\documentclass[a4paper,oneside,english]{smfcustom}
 \pdfoutput=1
\usepackage{amssymb,mathrsfs,mathtools}
\usepackage[english]{babel}
\usepackage[utf8]{inputenc}
\usepackage[T1]{fontenc}
\usepackage{etex,tikz}\usetikzlibrary{patterns}
\usepackage{xcolor}\usepackage[all]{xy}
\usepackage{enumitem}
\usepackage{xspace}
 \def\N{\mathbb{N}} 
\def\Zpos{\mathbb{Z}_{\scriptscriptstyle >0}}
\def\P{\mathbb{P}} % projectivization
\def\K{\mathbf{k}} % base field
\def\X{\mathcal{X}} % family of complete intersections
\def\Y{\mathcal{Y}} % universal Grassmannian
\def\G{\mathbf{G}} % parameter space of \Y
\def\Q{\mathcal{Q}} % Plücker line bundle on \G
\def\O{\mathcal{O}}
\def\B{\mathbf{B}} % base locus

\def\L{\mathcal{L}}
\def\pxi{[\xi]} % the line in the direction of xi
\let\leq\leqslant \let\geq\geqslant
\let\mathbi\boldsymbol
\def\inverse{^{\scriptscriptstyle-1}}
\def\mand{\qquad\text{and}\qquad}
\def\ie{\textit{i.e.}\xspace} \def\eg{\textit{e.g.}\xspace} \def\cf{\textit{cf.}\xspace}
\def\Wedge{\mathsf{\Lambda}}
\def\loceq{\stackrel{\text{loc}}=}
\let\bydef\coloneqq
\def\diff{\mathop{}\!{\mathrm{d}}}
\DeclareMathOperator{\Bs}{Bs} % Base locus
\DeclareMathOperator{\Grass}{Grass} % Grassmannian
\DeclareMathOperator{\pr}{pr} % projection
\DeclareMathOperator{\codim}{codim} % codimension
\DeclareMathOperator{\Span}{Span} % Span
\DeclareMathOperator{\rk}{rank} % rank
\DeclareMathOperator{\Exc}{Exc} % exceptional locus
\DeclareMathOperator{\nci}{\infty} % non finite locus
\DeclareMathOperator{\id}{id} %identity
\DeclareMathOperator{\res}{res} %restriction
\DeclarePairedDelimiter{\abs}{\lvert}{\rvert} % \abs* adjust size
\DeclarePairedDelimiter{\Set}{\{}{\}} % \Set* adjust size
\newcommand{\xym}[3][1]{\left.\vcenter{\xy\xymatrix"m"@R=.5pt@C=.5pt@W=1em@H=#1em{#3}\POS"m1,1"."m#2"!C*\frm{(}*\frm{)}\endxy}\right.}
\def\midbar{
  \mathchoice
  {\mathrel{\textsl{\Large|}}}
  {\mathrel{\textsl{\large|}}}
  {\mathrel{\textsl{|}}}
  {\mathrel{\textsl{\small|}}}
}
\def\displaymap#1#2{%
  \ifx\relax#1\relax
  \left.\vcenter{\xymatrix@=0pc{#2}}\right.
  \else
  #1\colon\left\vert\vcenter{\xymatrix@=0pc{#2}}\right.
  \fi
}

\def\POX{\P(\Omega_{X})}
\def\POM{\P(\Omega_{M})}
\def\POMI{\P(\Omega_{M}\rvert_{M_{I}})}
\def\POXS{\P(\Omega_{\X/S})}
\def\PkI{\P(\K_{I})}
\makeatletter
\def\a{\@ifnextchar^{\@a}{\@@a}}
\def\@a^#1{\mathbi{a}^{#1}}
\def\@@a{\mathbi{a}^{\bullet}}
\makeatother

\author{Damian Brotbek}
\address{Institut de Recherche Mathématique Avancée, Université de Strasbourg}
\email{brotbek@math.unistra.fr}
\author{Lionel Darondeau}
\address{Institut de Recherche Mathématique Avancée, Université de Strasbourg}
\curraddr{Currently: Departement Wiskunde, KU Leuven}
\email{lionel.darondeau@normalesup.org}
\thanks{The research of the second author was supported by the USIAS project \textit{``Rational Points, Rational Curves and Automorphisms of Special Varieties''} of Carlo Gasbarri and Gianluca Pacienza (http://www.usias.fr).}
\title[Complete intersection varieties with ample cotangent bundles]
{Complete intersection varieties\\with ample cotangent bundles}
\keywords{Ample cotangent bundle, symmetric differential forms, complete intersection varieties}
\subjclass{14M10}
\date{\today} %version actuelle
\begin{document}
%%%
%ENVS 
%%%
\newenvironment{THM}{\begin{enonce*}{Main result}}{\end{enonce*}}

\begin{abstract}
  Any smooth projective variety contains many complete intersection subvarieties with ample cotangent bundles, of each dimension up to half its own dimension.
\end{abstract}
\maketitle
It is expected that projective varieties with ample cotangent bundles should be abundant, at least under reasonable assumptions; in~\cite{Sch92}, Schneider proved that such a variety cannot be embedded in a projective space in such a way that its dimension is larger than its codimension.
Yet, not so many higher dimensional examples were known until recently.
Looking for such examples and generalizing a question of Schneider in the cited work, one can wonder if a smooth \(N\)-dimensional projective variety contains subvarieties with ample cotangent bundles (besides curves). 
In this paper, we answer this question and establish that, after taking into account Schneider's condition, subvarieties with ample cotangent bundles are ubiquitous.

\begin{THM}
  In every smooth projective variety \(M\), for each \(n\leq\dim(M)/2\), there exist smooth subvarieties of dimension \(n\) with ample cotangent bundles.
\end{THM}
To give some illustrative examples: in view of the aforementioned Schneider's result, this statement is sharp for \(M=\P^{N}\) (which has anti-ample cotangent bundle), 
by~\cite{Deb05} it is also sharp if \(M\) is an abelian variety (which has trivial cotangent bundle),
and the statement becomes obviously non-sharp (and trivial) if \(M\) itself has already ample cotangent bundle.

In an ambient variety \(M\), a natural way to construct subvarieties with ample cotangent bundles is to consider complete intersections of very ample divisors. 
Indeed, there are several ways to see that cotangent bundles of smooth hypersurfaces carry more ``positivity'' than the cotangent bundle of the original variety itself 
(consider \eg adjunction formula, or observe that cotangent bundles of hypersurfaces are quotients of the cotangent bundle of the ambient variety). Taking the complete intersection of more and more hypersurfaces, it is a natural question to ask at which point the cotangent bundle may become ample. 
We say that a property holds for a \textsl{general} member of some family, if 
it holds for any member of the family over a dense Zariski-open subset of the base.
In~\cite{Deb05}, Debarre conjectured that complete intersections in projective space having codimensions at least equal to their dimensions should have ample cotangent bundles, provided the hypersurfaces one intersects are general and sufficiently ample. This conjecture was motivated by his proof of the analogous statement for complete intersections in abelian varieties (see also~\cite{Deb05Cor,Ben11}). 

In the spirit of this conjecture, we prove the following statement, which implies our main result.
\begin{theo}
  \label{thm:ample}
  On a \(N\)-dimensional smooth projective variety \(M\), equipped with a very ample line bundle \(\O_{M}(1)\),
  if \(N/2\leq c\leq N\),
  for each \(\mathbi{\delta}=(\delta_{1},\dotsc,\delta_{c})\in(\Zpos)^{c}\),
  there exists \(\nu(\mathbi{\delta})\in\mathbb{Q}\)
  such that for all multi-degrees
  \(
  (d_{1},\dotsc,d_{c})
  =
  \nu\cdot(\delta_{1},\dotsc,\delta_{c})
  \)
  with \(\nu\in\mathbb{Q}\) such that
  \(
  \nu
  \geq
  \nu(\mathbi{\delta})
  \), 
  the complete intersection of general hypersurfaces 
  \(
  H_{1}\in\abs*{\O_{M}(d_{1})},\dotsc,H_{c}\in\abs*{\O_{M}(d_{c})}
  \)
  has ample cotangent bundle.
\end{theo}

Our proof of Theorem~\ref{thm:ample} has two inherent flaws, that are equally inconsequential regarding our main result. 
The first one is that we cannot use powers of different very ample divisors to define \(H_{1},\dotsc,H_{c}\). 
The second one is that the lower bound on the degrees is not uniform (namely it depends on the direction \(\mathbi{\delta}\)).
Moreover, our estimate of \(\nu(\mathbi{\delta})\) is only ``almost effective'', 
because it depends on some integer \(m(\mathbi{\delta})\) involved in some Noetherianity argument 
(see Theorem~\ref{thm:ampleStrong} in Sect.~\ref{sse:positivity}).

During the last steps of the preparation of this paper, we were informed that Xie announced a proof of the ampleness of cotangent bundles of general complete intersections in \(M=\P^{N}\) having codimensions at least as large as their dimensions, with a uniform lower bound on the degrees \(d_{1},\dotsc,d_{c}\gtrsim N^{N^{2}}\) (\cite{Xie15}). 
His result is thus stronger than Theorem~\ref{thm:ample}, concerning the hypothesis on the degrees.
Let us mention that his ``product coup'' could be adapted in our situation to give a uniform lower bound on the degrees, modulo a technical generalization of our arguments.

Soon after a preliminary version of the present paper was made available on arXiv, 
Deng was able to give an effective bound on \(m(\mathbi{\delta})\), which, in the particular case \(d_{1}=\dotsb=d_{c}=d\), allowed him to derive the effective lower bound \(d\geq 8c(2N)^{N+c}\).
Combining this with Xie's ``product coup'', Deng has made our proof fully effective, with the uniform lower bound 
\(
d_{1},\dotsc,d_{c}
\geq 
16c^{2}(2N)^{2N+2c}
\)
in Theorem~\ref{thm:ample} (see~\cite{Den16,Den17}).
\newline

Our result is in the vein of several previous ones we would like to mention. 
In an unpublished work, Bogomolov constructed complete intersection varieties with ample cotangent bundles in products of varieties with big cotangent bundles (we refer to~\cite{Deb05} for a treatment of this result). 
As already mentioned, Debarre proved in~\cite{Deb05} that in an abelian variety, a general complete intersection of sufficiently ample divisors, whose codimension is at least as large as its dimension, has ample cotangent bundle. 
In~\cite{Bro14}, the first author proved, among other results in the direction of the quoted conjecture of Debarre, that Theorem~\ref{thm:ample} (with an effective bound on the degrees) holds for complete intersection surfaces in projective space,
and in~\cite{Bro16} he also proved that a general complete intersection \(X\) of multi-degree \((d,\dotsc,d)\) in \(\P^{N}\), with \(\codim X\geq 3\dim X-2\), has ample cotangent bundle when \(d\geq 2N+3\).

Recall that the cotangent bundle \(\Omega_{X}\) of a manifold \(X\) is said \textsl{ample} if the canonical line bundle \(\O_{\POX}(1)\) on the projectivization \(\POX\) is itself ample%
\footnote{
  for a vector bundle \(E\) on a variety \(X\) we denote by \(\pi_{E}\colon\P(E)\to X\) the projectivization of rank one quotients of \(E\)
}.
The positivity of the cotangent bundle is thus related to the existence and the geometry of \textsl{symmetric differential forms}
\[
  \sigma
  \in
  H^{0}(\POX,\O_{\POX}(m))
  \cong
  H^{0}(X,S^{m}\Omega_{X}).
\]
In~\cite{Bro14}, Demailly's holomorphic Morse inequalities are used in order to prove that 
for a complete intersection \(X=H_{1}\cap\dotsb\cap H_{c}\subset M\) of sufficiently ample divisors, with \(\codim(X)\geq\dim(X)\),
\(\Omega_{X}\) is \textsl{big}. 
This means that there is a constant \(C>0\) such that for any \(m\) large enough,
\(
\dim
H^{0}(X,S^{m}\Omega_{X})
\geq 
C\cdot m^{\dim\POX}
\)
(see~\cite{Laz04II}).
However, the existence of many symmetric differential forms does not provide much information on their base locus.

In~\cite{Bro16}, a more computational approach was introduced in order to construct symmetric differential forms vanishing on an ample divisor,
with an explicit dependance in the equations of a complete intersection \(X\) in \(\P^{N}\).
One way to avoid the computational difficulties that  occur by following this approach
is to deduce the statement for general complete intersections from a well chosen particular case.
Indeed, ampleness being an open property in families, to prove that a general complete intersection of a given multi-degree is ample,
it is sufficient to show that there exists an example of a smooth  complete intersection of the required multi-degree having ample cotangent bundle. 
These considerations motivate the choice to work with a certain subfamily of the family of complete intersections---concretely: with complete intersections defined by homogeneous equations of a special form---.
In~\cite{Bro16}, the computation is performed for a subfamily of complete intersections defined by sufficiently many \textsl{deformations of Fermat type hypersurfaces} (see Sect.~\ref{sse:bihom} for a generalization).
This choice of equations was certainly motivated by the technical simplification it provides,
but also by the works~\cite{BG77,MN96,Nad89,DEG97,SZ02} and others on examples of Kobayashi hyperbolic varieties. 
Indeed, varieties with ample cotangent bundles are hyperbolic (\cite{Kob98}), 
hence it is natural to start with hypersurfaces that are likely to be hyperbolic. 

In his aforementionned work~\cite{Xie15}, Xie refines the paper~\cite{Bro16} by working with more elaborated deformations of Fermat type hypersurfaces.
A subtle choice of equations allows him to produce a greater number of tractable symmetric differential forms and to control their base locus.
The present work also finds its roots in~\cite{Bro16}, and pushes further the ideas initiated there to complete the study of families of complete intersection varieties.
In contrast with the cited works~\cite{Bro16,Xie15}, we do not manipulate explicit expressions of some symmetric differential forms.
Although our intuitions grew from the geometric interpretation of a generalization of the cohomological computations arising in~\cite{Bro16},
in the present text we choosed to emphasize the intrinsic simplicity and the geometric nature of the proof.

While to the best of our knowledge the presented approach (outlined in the next section) is new, it involves some general ideas that are already present in earlier works. 
The independence between the equations defining the complete intersection and their differentials is a key ingredient in Debarre's proof~\cite{Deb05}.
The use of the ``positivity'' of the parameter space of a family of varieties is present in many recent works on hyperbolic varieties subsequent to Voisin's article~\cite{Voi96}, although it is used there in a very different way.
\newline

Let us lastly mention that the first author has extended the techniques of the present work from the setting of symmetric differentials to the setting of higher order jet differentials with the aim of proving a hyperbolicity conjecture of Kobayashi (\cite{Bro17}).

\section{Outline of the proof of Theorem~\ref{thm:ample}}
\label{se:outline}
Our proof of Theorem~\ref{thm:ample} requires some technical considerations that might hide the underlying  geometry.  
For this reason, we wish to propose an outline of the proof, in which we first skip the technical details (to be provided in the remaining of the text).
It will hopefully contribute to the overall understanding of the arguments.

Assume for simplicity until the end of this section that the ambient variety \(M=\P^{N}\) is the projective space of dimension \(N\geq2\) over an algebraically closed field \(\K\) of characteristic \(0\).
Fix homogenous coordinates \([y_{0}:\dotso:y_{N}]\) in \(\P^{N}\).  
For \(J=(j_{0},\dotsc,j_{N})\in \N^{N+1}\), the standard multi-index notation \(\mathbi{y}^{J}\bydef y_{0}^{j_{0}}\dotsm y_{N}^{j_{N}}\) will be used.
For all unexplained notation used in this introductory section, concerning multi-indices \(J\in\N^{N+1}\),  \cf Sect.~\ref{se:ampleness}. 

\subsection{Choice of the equations of \ensuremath{\X} and local equations of \ensuremath{\POXS}}
Consider complete intersections 
\(
X
\bydef 
H_{1}\cap\dotsb\cap H_{c}
\)
of codimension \(c\),
where each hypersurface \(H_{p}\) (\(p=1,\dotsc,c\)) is defined, 
for suitable integers \(\varepsilon_{p}>0\), \(\delta_{p}>0\) and 
\(r\in \N\),
by a homogenous polynomial of the form
\begin{equation}
  \label{eq:shape}
  E_{p}(\a^{p})
  \bydef
  \sum_{\abs{J}=\delta_{p}}
  a_{J}^{p}\,
  \mathbi{y}^{(r+1)J},
\end{equation}
with
\(a_{J}^{p}\in H^{0}\big(\P^{N},\O_{\P^{N}}(\varepsilon_{p})\big)\) and where we use the vectorial notation \(\a^{p}\bydef(a_{J}^{p})_{\abs{J}=\delta_{p}}\). 
Accordingly, 
\(
E_{p}(\a^{p})
\in 
H^{0}\big(\P^{N},\O_{\P^{N}}(\varepsilon_{p}+(r+1)\delta_{p})\big)
\)
(see Sect.~\ref{sse:bihom}).
This choice of equations is motivated by the previous work~\cite{Bro16} by the first author, where the case \(\delta_{1}=\dotsb=\delta_{c}=1\) was studied.

For \(\delta\in\N\), denote by 
\[
  N_{\delta}
  \bydef
  \#\Set{\abs{J}=\delta}
  = 
  \dim H^{0}\big(\P^{N},\O_{\P^{N}}(\delta)\big)
\] 
the number of parameters of a degree \(\delta\) homogenous equation in \(N+1\) variables.
Let 
\[
  S
  \subset
  H^{0}\big(\P^{N},\O_{\P^{N}}(\varepsilon_{1})\big)^{N_{\delta_{1}}}
  \times
  \dotsb
  \times
  H^{0}\big(\P^{N},\O_{\P^{N}}(\varepsilon_{c})\big)^{N_{\delta_{c}}}
\]
be the subspace of parameters \((\a^{1},\dotsc,\a^{c})\) of equations \((E_{1}(\a^{1}),\dotsc,E_{c}(\a^{c}))\) as above defining smooth complete intersections.
Denote by \(\X\to S\) the universal family of such complete intersections, defined in \(S\times\P^{N}\) by the universal equations \((E_{1},\dotsc,E_{c})\), and consider the projectivization of the relative cotangent bundle:
\[
  \pi_{\Omega_{\X/S}}
  \colon
  \POXS
  \longrightarrow
  \X.
\]
Local equations for \(\POXS\subset S\times\P(\Omega_{\P^{N}})\) are given by the universal equations \((E_{1},\dotsc,E_{c})\) and their relative differentials \((\diff_{\P^{N}}E_{1},\dotsc,\diff_{\P^{N}}E_{c})\).

From our point of view, the key feature of equations~\eqref{eq:shape} is that such equations and their differentials can be written in the same formal shape. 
Locally on \(\P^{N}\), for \(p=1,\dotsc,c\), one can write
\begin{equation}
  \label{eq:cotangent}
  E_{p}
  =
  \sum_{\abs{J}=\delta_{p}}
  \alpha_{J}^{p}\,
  (\mathbi{y}^{r})^{J}
  \mand
  \diff_{\P^{N}} E_{p}
  \loceq
  \sum_{\abs{J}=\delta_{p}}
  \theta_{J}^{p}\,
  (\mathbi{y}^{r})^{J},
\end{equation}
where the notation \(\loceq\) emphasizes that this makes sense only locally,
after restricting ourselves to some affine open subset \(V\subset\P^{N}\) and trivializing \(\O_{\P^{N}}(1)\rvert_{V}\).
Here the ``coefficients'' \(\alpha_{J}^{p}\) are homogenous polynomials of degree \(\varepsilon_{p}+\delta_{p}\),
and the ``coefficients''
\(\theta_{J}^{p}\) are \(1\)-forms on \(V\) depending on our choice of trivialization for \(\O_{\P^{N}}(1)\rvert_{V}\) (see Sect.~\ref{sse:reduction}). 

\subsection{The morphism \(\varPsi\)}
After making the substitution \(z_{i}=y_{i}^{r}\) for \(0\leq i\leq N\), 
the equations in~\eqref{eq:cotangent} become
\[
  \sum_{\abs{J}=\delta_{p}}
  \alpha_{J}^{p}\,\mathbi{z}^{J}
  \mand
  \sum_{\abs{J}=\delta_{p}}
  \theta_{J}^{p}\,\mathbi{z}^{J}
  \qquad{(p=1,\dotsc,c)}.
\]
Thinking this time of \(\alpha_{J}^{p}\) and \(\theta_{J}^{p}\) as variables in \(\K\) provide us with a way to compare \(\POXS\) with the universal family \(\Y\) of ``complete intersections'' of multi-degree \((\delta_{1},\delta_{1},\dotsc,\delta_{c},\delta_{c})\) in \(\P^{N}\) defined by these universal equations.
We chose to parametrize \(\Y\) not by
\((\prod_{p=1}^c\K^{N_{\delta_p}}\times \K^{N_{\delta_p}})\)
but rather by a product of Grassmannians
\[
  \G
  \bydef
  \Grass\big(2,\K^{N_{\delta_{1}}}\big)
  \times
  \dotsb
  \times
  \Grass\big(2,\K^{N_{\delta_{c}}}\big),
\]
and we therefore have  \(\Y\subset\G\times\P^{N}\) (see~\eqref{eq:Y}). 
The choice of this parameter space will be justified shortly.

In Section~\ref{sse:reduction} we precise the relation between \(\POXS\) and \(\Y\) by defining, 
for some suitable open subset \(S^{\circ}\subset S\), a morphism \(\varPsi\) 
from the projectivized relative cotangent bundle \(\P(\Omega_{\X/S^{\circ}})\) 
to \(\Y\):
\[
  \xymatrix{
    \P(\Omega_{\X/S^{\circ}})\ar[d]\ar[rr]^{\varPsi}&&\Y\ar[d]_{\rho}\ar[r]&\P^{N}.
    \\
    S^{\circ}&&\G&
  }
\]
The above diagram is detailed in~\eqref{eq:Diagram}.
To construct \(\varPsi\), it is critical to parametrize \(\Y\) by a product of Grassmannians rather than by a product of projective spaces, because we can only define \(\diff_{\P^{N}} E_{p}\) locally (as explained in Sect.~\ref{sse:reduction}). 
The restriction to the open set \(S^{\circ}\) is needed to work with sufficiently independant parameters \(\alpha_{J}^{p}\) and \(\theta_{J}^{p}\) (see Prop.~\ref{prop:openS})
and the morphism \(\varPsi\) is then defined by~\eqref{eq:DefinitionPsi}.

\subsection{The model situation}
The product of Grassmannians \(\G\) is naturally endowed with an ample line bundle \(Q\to\G\), 
namely the tensor product of the Pl\"{u}cker line bundles coming from the different factors%
\footnote{
  Given varieties \(X_{1}\) and \(X_{2}\), vector bundles \(E_{1}\) on \(X_{1}\) and \(E_{2}\) on \(X_{2}\), 
  one can define the tensor product \(E_{1}\boxtimes E_{2}\bydef p_{1}^*E_{1}\otimes p_{2}^*E_{2}\) where \(p_{i}\colon X_{1}\times X_{2}\to X_{i}\), \(i=1,2\), are the canonical projections. 
}.
Denote \(\O_{\X}(1)=\big(\O_{S}\boxtimes\O_{\P^{N}}(1)\big)\rvert_{\X}\).
By the very construction of the morphism \(\varPsi\), the pullback of
\((Q^{\otimes m}\boxtimes\O_{\P^{N}}(-1))\rvert_{\Y}\) is of the form:
\begin{equation}
  \label{eq:pullback}
  \varPsi^{\ast}\big((Q^{\otimes m}\boxtimes\O_{\P^{N}}(-1))\rvert_{\Y}\big)
  =
  \O_{\P(\Omega_{\X/S^{\circ}})}(cm)
  \otimes 
  \pi_{\Omega_{\X/S^{\circ}}}^{\ast}\O_{\X}(\star-r),
\end{equation}
for some integer \(\star\) (computed later) \emph{independent of \(r\)} (see~\eqref{eq:PullBack}).
Hence, for \(r>\star\), and for \(s\in S^{\circ}\), every global section of \((Q^{\otimes m}\boxtimes\O_{\P^{N}}(-1))\rvert_{\Y}\) gives rise to a symmetric differential form on \(X_{s}\) vanishing along some ample divisor. 

Suppose now \(2c\geq N\) like in the assumptions of Theorem~\ref{thm:ample}.
The morphism \(\rho\colon\Y\to\G\) is then generically finite. 
The pullback \(\rho^{\ast}Q\) of the ample line bundle \(Q\) is not ample on \(\Y\), but it is still big and nef.
Thus Nakamaye's Theorem on the augmented base locus (\cite{Nak00}) provides a \emph{geometric control on its base locus}.
Namely, for all sufficiently large integers \(m\) one has
\begin{equation}
  \label{eq:NakaNaka}
  \Bs\big(
  (Q^{\otimes m}\boxtimes\O_{\P^{N}}(-1))\rvert_{\Y}
  \big)
  =
  \Exc(\Y\to\G)
  \subset 
  \mathcal{Y},
\end{equation}
where \(\Exc(\Y\to\G)\) is  the reunion of all positive dimensional fibers of the morphism \(\mathcal{Y}\to\G\).

Building on some ideas of the work of Benoist~\cite{Ben11} (see Sect.~\ref{se:Olivier}) one shows that the locus \(\Exc(\Y\to \G)\) is small (\ie of large codimension) when \(\delta_{1},\dotsc,\delta_{c}\) are large. 

\subsection{Pulling back the positivity from \ensuremath{\Y} to \ensuremath{\P(\Omega_{\X/S^{\circ}})}}
Hence, for \(m\) satisfying~\eqref{eq:NakaNaka} and for \(s\in S^{\circ}\),
the combination of~\eqref{eq:pullback} and~\eqref{eq:NakaNaka} 
gives us some information on the base locus of the line bundle 
\[
  \L_{s}
  \bydef
  \O_{\P(\Omega_{X_{s}})}(cm)\otimes \pi_{\Omega_{X_{s}}}^*\O_{X_{s}}(\star-r).
\]
Our goal is then to prove that for such \(m\), for \(r>\star\),
and for a general element \(s\in S\), \(\L_{s}\)
is nef. 
It would indeed follow that \(\Omega_{X_{s}}\) is ample for a general member \(X_{s}\) of the family \(\X\) of special complete intersections under consideration.
By the openness property of ampleness, this would end the proof of Theorem~\ref{thm:ample} for suitable multi-degrees (see Sect.~\ref{sse:condition}).

Observe that the range of integers \(m\) satisfying the condition~\eqref{eq:NakaNaka} depends only on the geometry of \(\Y=\Y(\mathbi{\delta})\), which is clearly independent of \(r\).
Since the integer \(\star\) depends in turn on \(m\), it is in the contrary the range of adequate parameters \(r>\star\) that will depend on the geometry of \(\Y\), and we can safely fix \(r\) \textit{a posteriori}.

Our idea for obtaining the nefness of \(\L_{s}\) is to  control the inverse image under \(\varPsi\) of the locus \(\Exc(\Y\to \G)\);
If for general \(s\) one could prove that for any curve \(\mathscr{C}\subset \P(\Omega_{X_{s}})\), one has
\[
  \varPsi(\mathscr{C})
  \not\subset
  \Exc(\Y\to\G)
  =
  \Bs\big((Q^{\otimes m}\boxtimes\O_{\P^{N}}(-1))\rvert_{\Y}\big),
\]
then, in view of~\eqref{eq:pullback} and~\eqref{eq:NakaNaka}, it would follow that for any curve \(\mathscr{C}\subset \P(\Omega_{X_{s}})\), one has \(c_{1}(\L_{s})\cdot \mathscr{C}\geq 0\), as desired.

Regretably, the particular shape of our equations causes some technical complications along the coordinate hyperplanes in \(\P^{N}\).
This prevents us to directly apply our idea and forces us to consider the stratification on \(\P(\Omega_{\P^{N}})\) induced by these coordinate hyperplanes in \(\P^{N}\).
Then, studying each stratum independently, by modifying the above approach in some minor ways, and putting the information of the different strata together allows us to prove the sought nefness.
See Sect.~\ref{sse:exc},~\ref{sse:positivity} for full detail. 

\section{Ampleness of the cotangent bundle of general complete intersections}
\label{se:ampleness}
We work over an algebraically closed field \(\K\) of characteristic \(0\).  Recall that \(M\) is an \(N\)-dimensional smooth projective variety over \(\K\), equipped with a very ample line bundle \(\O_{M}(1)\). We always assume \(N\geq2\).

Throughout this text we use the following notation.
We denote by \(\Set{\abs{J}=\delta}\subset\N^{N+1}\) the subset of multi-integers \(J=(j_{0},\dotsc,j_{N})\) such that
\(\abs{J}\bydef j_{0}+j_{1}+\dotsb+j_{N}=\delta\).
There are
\(
N_{\delta}
\bydef
\binom{N+\delta}{N}
\)
such multi-indices.
The support of a multi-index \(J\) is denoted
\[
  [J]
  \bydef
  \Set*{
    i\in\Set{0,1,\dotsc,N}\midbar j_{i}\neq0
  }.
\]
For \(I\subsetneq\Set{0,\dotsc,N}\) with \(\#I=(N-k)\), there are \(\binom{k+\delta}{k}\) multi-indices of length \(\delta\) with \([J]\cap I=\varnothing\).
These parametrize degree \(\delta\) monomials that do not vanish identically on the vector space
\[
  \K_{I}
  \bydef
  \Set*{(z_{0},\dotsc,z_{N})\in\K^{N+1}\midbar z_{i}=0,\forall i\in I}.
\]
Considering the induced homogenous coordinates \([z_{0}:\dotso:z_{N}]\) on \(\P^{N}\cong\P(\K^{N+1})\), one has isomorphisms
\begin{equation}
  \label{eq:Identification}
  H^{0}(\P^{N},\O_{\P^{N}}(\delta))
  \cong 
  \K^{\binom{N+\delta}{N}}
  \cong 
  \bigoplus_{\abs{J}=\delta}\K
  \mand
  H^{0}(\PkI,\O_{\PkI}(\delta))
  \cong 
  \K^{\binom{k+\delta}{k}}
  \cong 
  \bigoplus_{\mathclap{\abs{J}=\delta\colon[J]\cap I=\varnothing}}\K.
\end{equation}

When we consider points in \(\POM\), we will use the following convention: when we write \(\pxi\in\POM\) we tacitly assume, for \(x=\pi_{\Omega_{M}}(\pxi)\), that we have fixed \(\xi\in T_{x}M\setminus\Set{0}\) representing \(\pxi\). 
Moreover, given a function \(f\) in a neighborhood of \(x\) we will write \(f(\xi)\) instead of \(f(x)\), only implying the pre-composition by the projection \(\xi\mapsto x\).

When we work on a chart \(V\subset M\) in a trivializing affine covering \(\mathfrak{V}\) for \(\O_{M}(1)\), given \(\sigma\in H^{0}(M,\O_{M}(1))\) we will denote by \((\sigma)_{V}\in \O(V)\) the regular function corresponding to \(\sigma\) under that trivialization.

\subsection{Complete intersections cut out by bihomogeneous sections}
\label{sse:bihom}
We fix \(N+1\) sections in general position,
\(
\zeta_{0},\dotsc,\zeta_{N}
\in 
H^{0}\big(M,\O_{M}(1)\big),
\)
and we set \(D_{i}\bydef(\zeta_{i}=0)\) for \(i=0,\dotsc,N\). 
By ``general position'', we mean that each of the \(D_{i}\) is smooth, and that the divisor \(D=\sum_{i}D_{i}\) is simple normal crossing.
In the familiar case \(M=\P^{N}\), this means that \(\zeta_{0},\dotsc,\zeta_{N}\) are homogenous coordinates in \(\P^{N}\).

For an integer \(c\) satisfying the hypothesis \(N/2\leq c\leq N\) of Theorem~\ref{thm:ample}, 
and two \(c\)-tuples of positive integers
\(\mathbi{\varepsilon}=(\varepsilon_{1},\dotsc,\varepsilon_{c})\),
\(\mathbi{\delta}=(\delta_{1},\dotsc,\delta_{c})\),
we construct the family \(\X\) mentioned in Sect.~\ref{se:outline} as follows.
For \(p=1,\dotsc,c\), 
for
\[
  \a^{p}
  =
  \big(
  a_{J}^{p}
  \in
  H^{0}\big(M,\O_{M}(\varepsilon_{p})\big)
  \big)_{\abs{J}=\delta_{p}},
\]
and for a positive integer \(r\) fixed later according to our needs, 
we consider the bihomogeneous section of \(\O_{M}(\varepsilon_{p}+(r+1)\delta_{p})\) over \(M\) defined by\footnote{here and throughout the text, we use the standard multi-index notation for multivariate monomials}
\[
  E^{p}(\a^{p},\cdot)
  \colon
  x
  \mapsto
  \sum_{\abs{J}=\delta_{p}}
  a_{J}^{p}(x)
  \mathbi{\zeta}(x)^{(r+1)J}.
\]
We rather want to let \(\a^{p}\) vary in the parameter space
\[
  S_{p}
  \bydef 
  H^{0}\big(M,\O_{M}(\varepsilon_{p})\big)^{N_{\delta_{p}}}
\]
and to think at \(E^{p}\) as a section
\(E^{p}\in H^{0}(S_{p}\times M,\O_{S_{p}}\boxtimes\O_{M}(d_{p}))\)   
where \( d_{p} \bydef \varepsilon_{p}+(r+1)\delta_{p}\).

We then consider the universal family \(\overline{\X}\subset S_{1}\times\dotsb\times S_{p}\times M\) defined by those sections, \ie
\[
  \overline{\X}
  \bydef
  \Set*{
    (\a^{1},\dotsc,\a^{c};x)
    \in
    S_{1}\times\dotsb\times S_{c}\times M
    \midbar
    E^{1}(\a^{1},x)=0,\dotsc, E^{c}(\a^{c},x)=0
  }.
\]

\begin{lemm}
  The general fibers of the family \(\overline{\X}\to S_{1}\times\dotsb\times S_{c}\) are smooth.
\end{lemm}
\begin{proof}
  It is sufficient to produce one smooth member of this family. To do so, it is possible to construct the sought complete intersection step-by-step, using the following statement (established below): 
  \begin{quote}
    \itshape
    Given \(\varepsilon,\delta >0\) and a smooth subvariety \(X\) of \(M\)  there exists a global section \(\sigma\) of \(\O_{M}(\varepsilon +(r+1)\delta)\) of the form
    \(\sigma =\sum_{i=0}^{N}a_{i}\zeta_{i}^{(r+1)\delta}\), where \(a_{0},\dotsc,a_{N}\in H^{0}(M,\O_{M}(\varepsilon))\), such that \(X\cap (\sigma=0)\) is smooth.
  \end{quote}
  The proof is a minor modification of the proof of the classical Bertini Theorem, given here for the sake of completeness.

  Consider the set \(\Sigma\subset H^{0}(M,\O_{M}(\varepsilon))^{N+1}\times X\) composed of element of the form \((a_{0},\dotsc,a_{N},x)\) such that the hypersurface \(H_{\mathbi{a}}\) defined by the vanishing of \(\sigma_{\mathbi{a}} \bydef\sum_{i=0}^{N}a_{i}\zeta_{i}^{(r+1)\delta}\) contains \(x\) but \(X\cap H_{\mathbi{a}}\) has a singularity at \(x\). 
  Denote by \(p_{2}\colon\Sigma\to X\) the canonical projection.
  We are going to estimate, for all \(x\in X\), the dimension of the fibers \(\Sigma_{x}\bydef p_{2}\inverse(x)\).

  Fix \(x\in X\), set \(H^{0}(M,\O_{M}(\varepsilon))_{x}^{N+1}=\Set{\mathbi{a}\midbar x\in H_{\mathbi{a}}}\), this is a hyperplane of \(H^{0}(M,\O_{M}(\varepsilon))^{N+1}\).
  Fix also an affine neighborhood \(V\subset M\) of \(x\), on which \(\O_{M}(1)\) is trivialized. Then we have a linear map 
  \[
    \gamma_{x}
    \colon
    H^{0}(M,\O_{M}(\varepsilon))_{x}^{N+1}
    \longrightarrow
    \Omega_{X\cap V,x}, 
  \]
  defined by \(\gamma_{x}(\mathbi{a})=\diff(\sigma_{\mathbi{a}})_{V}\rvert_{\Omega_{X,x}}\). 
  Observe that \(\diff(\sigma_{\mathbi{a}})_{V}\in \Gamma(V,\Omega_{V})\), so it makes sense to restrict it to \( \Omega_{X\cap V,x}\). Now \(\Sigma_{x}\cong\ker\gamma_{x}\). Let us prove that \(\gamma_{x}\) is surjective, by showing that for \(\xi \in T_{x}X\setminus\Set{0}\), there exists \(\mathbi{a}\) such that 
  \(\sigma_{\mathbi{a}}(x)=0\) and
  \(\gamma_{x}(\mathbi{a})(\xi)\neq 0\). 
  There exists at least one index \(i_{0}\) such that \(\zeta_{i_{0}}(x)\neq0\).
  Since \(\varepsilon>0\), the line bundle \(\O_{M}(\varepsilon)\) is very ample and its sections separate first order jets.
  Hence, for each \(\xi\in T_{x}X\setminus\Set{0}\), 
  we can take \(\mathbi{a}\) such that \(a_{i}(x)=0\) for all indices \(i\), \(\diff(a_{i_{0}})_{V}\rvert_{x}(\xi)\neq0\) and \(\diff(a_{i})_{V}\rvert_{x}(\xi)=0\) for \(i\neq i_{0}\). A direct computation then shows that \(\gamma_{x}(\mathbi{a})(\xi)=d(a_{i_{0}})_{V}\rvert_{x}(\xi)\zeta_{i_{0}}^{(r+1)\delta}(x)\neq 0\), and this proves the desired surjectivity. 
  Therefore, 
  \[
    \dim \Sigma_{x}
    =
    \dim H^{0}(M,\O_{M}(\varepsilon))^{N+1}-1-\rk\gamma_{x}
    =
    \dim H^{0}(M,\O_{M}(\varepsilon))^{N+1}-1-\dim X.
  \]
  We infer that \(\dim(\Sigma)<\dim H^{0}(M,\O_{M}(\varepsilon))^{N+1}\).
  This implies that \(\Sigma\) does not dominate \(H^{0}(M,\O_{M}(\varepsilon))^{N+1}\), whence the result.
\end{proof}
We will from now on restrict ourselves  to  the dense open subset  \(S\subset S_{1}\times\dotsb\times S_{c}\) parametrizing smooth complete intersection varieties,
and we will prove the ampleness of the cotangent bundle of the general members of the universal family
\[
  \X
  \bydef
  \Set*{
    (\a^{1},\dotsc,\a^{c};x)
    \in
    S\times M
    \midbar
    E^{1}(\a^{1},x)=0,\dotsc, E^{c}(\a^{c},x)=0
  }.
\]
We denote by \(\pr_{1}\) and by \(\pr_{2}\) the natural projections from \(S\times M\) to its factors, 
and we denote by \(\pr_{1,\X}\) and by \(\pr_{2,\X}\) their restrictions to \(\X\).
Lastly, we denote
\(
\O_{\X}(1)
\bydef
\pr_{2,\X}^{\ast}\O_{M}(1)
\), 
for later use.

After projectivization, the relative cotangent sheaf \(\Omega_{\X/S}\) of the family \(\X\subset S\times M\) gives rise to a family
\[
  \xymatrix@C=0pt{
    \POXS
    \ar[d]_{\pi_{\Omega_{\X/S}}}&
    \subset &
    S\times\POM
    \ar[d]
    \\
    \X&
    \subset&
    S\times M.
  }
\]
We denote by \(\pr_{1}^{[1]}\) and by \(\pr_{2}^{[1]}\) the natural projections from \(S\times\POM\) to its factors,  
and we denote by \(\pr^{[1]}_{1,\X}\) and by \(\pr^{[1]}_{2,\X}\) their restrictions to \(\POXS\).

Since \(\pr_{1,\X}^{[1]}={\pr_{1,\X}}\circ{\pi_{\Omega_{\X/S}}}\), the fiber of \(\a\in S\) is \((\pr_{1,\X}^{[1]})\inverse(\a)=\P(\Omega_{X_{\a}})\) 
where \(X_{\a}\bydef(\pr_{1,\X})\inverse(\a)\).
Recall again that the cotangent bundle of \(X_{\a}\) is said ample if the canonical line bundle \(\O_{\P(\Omega_{X_{\a}})}(1)\) on \(\P(\Omega_{X_{\a}})\) is itself ample.

\subsection{Construction of families of non-positive-dimensional subschemes}
\label{sse:reduction}
We work temporarily on a chart \(V\) in a trivializing affine covering \(\mathfrak{V}\) for \(\O_{M}(1)\). 
The family \(\POXS\subset S\times\POM\) is  locally defined by the equations 
\(
E^{1}, \diff_{M}E^{1}, \dotsc, E^{c}, \diff_{M} E^{c}
\). 

The main geometric idea in our proof is the construction of a family \(\Y\) of non-positive-dimensional subschemes of \(\P^{N}\) defined by universal equations,
together with a map \(\varPsi\colon\POXS\to\Y\) used to ``pullback positivity''.
We get a number of equations that outreaches the dimension if we use both the equations \(E^{p}\) and their (relative) differentials \(\diff_{M} E^{p}\)---recall that \(2c\geq N\)---.

We think, for fixed \(\a^{p}\in S_{p}\), at \(E^{p}(\a^{p},\cdot)\) and \(\diff E^{p}(\a^{p},\cdot)\) as two independent polynomials of degree \(\delta_{p}\) in the variables \(\zeta_{0}^{r},\dotsc,\zeta_{N}^{r}\).
Namely, for \(\a^{p}\in S_{p}\) we consider locally the equation
\begin{align*}
  (E^{p}(\a^{p},\cdot))_{V}
  &=
  \sum_{\abs{J}=\delta_{p}}
  (\alpha_{J}^{p}(\a^{p},\cdot))_{V}
  \big(\mathbi{\zeta}\big)_{V}^{rJ}
  \in\O(V)
  \\\intertext{and its differential}
  \diff(E^{p}(\a^{p},\cdot))_{V}
  &=
  \sum_{\abs{J}=\delta_{p}}
  (\theta_{J}^{p}(\a^{p},\cdot))_{V}
  \big(\mathbi{\zeta}\big)_{V}^{rJ}
  \in H^{0}(V,\Omega_{V}),
\end{align*}%
where
\begin{equation}
  \label{eq:theta}
  \left\{
    \begin{aligned}
      \alpha_{J}^{p}(\a^{p},\cdot)\phantom{)_{V}}
      &\bydef
      a_{J}^{p}\mathbi{\zeta}^{J}
      &&\in
      H^{0}(M,\O_{M}(\varepsilon_{p}+\delta_{p}));
      \\
      (\theta_{J}^{p}(\a^{p},\cdot))_{V}
      &\bydef
      \diff(a_{J}^{p})_{V}
      (\mathbi{\zeta}^{J})_{V}
      +
      (a_{J}^{p})_{V}(r+1)\diff(\mathbi{\zeta}^{J})_{V} 
      &&\in 
      H^{0}(V,\Omega_{V}).
    \end{aligned}
  \right.%}
\end{equation}
The local sections \(\theta_{J}^{p}\) will soon disappear in favor of global sections that these will help to construct.

For \(\pxi\in\pi_{\Omega_{M}}\inverse(V)\) and for \(\a^{p}\in S_{p}\) set
\[
  (\mathbi{\alpha}^{p})_{V}(\a^{p},\xi)
  \bydef
  \big(\alpha_{J}^{p}(\a^{p},\xi)\big)_{\abs{J}=\delta_{p}}
  \in 
  \K^{N_{\delta_{p}}}
  \mand
  {(\mathbi{\theta}^{p})}_{V}(\a^{p},\xi)
  \bydef
  \big(\theta_{J}^{p}(\a^{p},\xi)\big)_{\abs{J}=\delta_{p}}
  \in 
  \K^{N_{\delta_{p}}}.
\]
As we shall see shortly, although the two points 
\((\mathbi{\alpha}^{p})_{V}(\a^{p},\xi)\)
and
\({(\mathbi{\theta}^{p})}_{V}(\a^{p},\xi)\)
in \(\K^{N_{\delta_{p}}}\) 
depend on the choice of trivialization and of the choice of a representative \(\xi\) for \(\pxi\), the vector space \(\Delta^{p}(\a^{p},\pxi)\) spanned by them does not. 
Under the assumption that this space is two dimensional, we get a well defined element in the Grassmannian of \(2\)-dimensional linear subspaces of \(\K^{N_{\delta_{p}}}\)
\[
  \Grass(2,N_{\delta_{p}})
  \bydef
  \Grass\big(2,\K^{N_{\delta_{p}}}\big)
  \cong 
  \Grass\big(2,H^{0}(\P^{N},\O_{\P^{N}}(\delta_{p}))\big).
\] 

In order to prove this, let us first recall the following classical construction of Wrońskian differential forms. 
\begin{lemm}
  Let \(L\to X\) be a line bundle over a smooth variety \(X\), and let \(\sigma_{1},\sigma_{2}\in H^{0}(X,L)\) be global sections. 
  The Wrońskian section
  \[
    \omega(\sigma_{1},\sigma_{2}) 
    \bydef 
    \det
    \xym{2,2}{\sigma_{1}&\sigma_{2}\\
    \diff\sigma_{1}&\diff\sigma_{2}}
  \]
  defines a global section in 
  \(H^{0}(X,\Omega_{X}\otimes L^{\otimes 2})
  \cong 
  H^{0}\big(\POX,\O_{\POX}(1)\otimes\pi_{\Omega_{X}}^{\ast}L^{\otimes 2}\big)\).
\end{lemm}
\begin{proof}
  The quotient \(\sigma_{1}/\sigma_{2}\) being a rational function on \(X\), 
  we can consider its logarithmic differential, which is a logarithmic form with simple poles along the divisor \((\sigma_{1}\sigma_{2}=0)\). Multiplying by \(\sigma_{1}\sigma_{2}\) cancels out the poles, and yields a section
  in \(H^{0}(X,\Omega_{X}\otimes L^{2})\).
  A short computation shows that this is precisely the Wrońskian form of the statement.
\end{proof}

Let us now fix \(1\leq p\leq c\) and an element \(\a^{p}\in S_{p}\), 
an open set \(V\subset M\) as above and \(\pxi\in\pi_{\Omega_{M}}\inverse(V)\).
Consider 
for \(J_{1},J_{2}\in\Set{\abs{J}=\delta_{p}}\)
the coordinates of
\(
(\mathbi{\alpha}^{p})_{V}(\a^{p},\xi)
\wedge
{(\mathbi{\theta}^{p})}_{V}(\a^{p},\xi)
\)
in
\(\Wedge^{2}\K^{N_{\delta_{p}}}\):
\begin{equation}
  \label{eq:defDeltaJJ}
  \Delta_{J_{1},J_{2}}^{p}(\a^{p},\xi)
  \bydef
  \det
  \xym{2,2}{
    {(\alpha_{J_{1}}^{p}(\a^{p},\xi))}_{V}&
    {(\alpha_{J_{2}}^{p}(\a^{p},\xi))}_{V}\\
    {(\theta_{J_{1}}^{p}(\a^{p},\xi))}_{V}&
    {(\theta_{J_{2}}^{p}(\a^{p},\xi))}_{V}
  }.
\end{equation}
\textit{A priori}, \(\Delta_{J_{1},J_{2}}^{p}(\a^{p},\cdot)\in \Gamma(V,\Omega_{V})\cong\Gamma(\cramped{\pi_{\Omega_{M}}\inverse(V)},\O_{\POM}(1))\) is only a local section of \(\O_{\POM}(1)\), but by an elementary linear algebra computation on \(\pi_{\Omega_{M}}\inverse(V)\) (recall~\ref{eq:theta}):
\[
  \omega\big(
  a_{J_{1}}^{p} \mathbi{\zeta}^{(r+1)J_{1}},
  a_{J_{2}}^{p} \mathbi{\zeta}^{(r+1)J_{2}}
  \big)
  \loceq
  \mathbi{\zeta}^{r(J_{1}+J_{2})}
  \Delta_{J_{1},J_{2}}^{p}(\a^{p},\cdot).
\]
Since
\[
  \omega\big(
  a_{J_{1}}^{p} \mathbi{\zeta}^{(r+1)J_{1}},
  a_{J_{2}}^{p} \mathbi{\zeta}^{(r+1)J_{2}}
  \big)
  \in 
  H^{0}\big(\POM,\O_{\POM}(1)\otimes\pi_{\Omega_{M}}^{\ast}\O_{M}(2 d_{p})\big).
\]
and 
\[
  \mathbi{\zeta}^{r(J_{1}+J_{2})}
  \in
  H^{0}\big(\POM,\pi_{\Omega_{M}}^{\ast}\O_{M}(2r\delta_{p})\big)
\]
are global sections,
one infers that actually
\[
  \Delta_{J_{1},J_{2}}^{p}(\a^{p},\cdot)
  \in
  H^{0}\big(\POM,\O_{\POM}(1)\otimes\pi_{\Omega_{M}}^{\ast}\O_{M}(2\varepsilon_{p}+2\delta_{p})\big).
\]
Moreover, \(\Delta_{J_{1},J_{2}}^{p}(\a^{p},\cdot)\) varies algebraically with respect to \(\a^{p}\),
it can thus be viewed as a section
\begin{equation}
  \label{eq:DeltaJJ}
  \Delta_{J_{1},J_{2}}^{p}
  \in 
  H^{0}\big(
  S\times\POM,
  (\pr_{2}^{[1]})^{\ast}\big(\O_{\POM}(1)\otimes\pi_{\Omega_{M}}^{\ast}\O_{M}(2\varepsilon_{p}+2\delta_{p})\big)
  \big).
\end{equation}

Passing to the projectivization \(\P(\Wedge^{2}\K^{N_{\delta_{p}}})\), we get a rational map  
\(
\Delta^{p}
\colon
S_{p}\times\POM
\dashrightarrow
\P(\Wedge^{2}\K^{N_{\delta_{p}}})
\)
defined in homogenous coordinates by
\(
\Delta^{p}
(\a^{p},\pxi)
=
\Big[
  \Delta_{J_{1},J_{2}}^{p}(\a^{p},\pxi)\midbar J_{1},J_{2}\in\Set{\abs{J}=\delta_{p}}
\Big].
\)
The indeterminacy locus of \(\Delta^{p}\) is contained in the base locus of the family \(\Set*{\Delta_{J_{1},J_{2}}^{p}\midbar{J_{1},J_{2}\in\Set{\abs{J}=\delta_{p}}}}\) which we will denote by \(B_{p}\).

Observe that by construction~\eqref{eq:defDeltaJJ}, outside \(B_{p}\), with the above notation,
\begin{equation}
  \label{eq:defDelta}
  \Delta^{p}(\a^{p},\pxi)
  \loceq
  \Span
  \big(
  {(\mathbi{\alpha}^{p})}_{V}(\a^{p},\xi),
  {(\mathbi{\theta}^{p})}_{V}(\a^{p},\xi)
  \big)
  \in 
  \Grass(2,N_{\delta_{p}})
  \subset 
  \P(\Wedge^{2} \K^{N_{\delta_{p}}}).
\end{equation}
Hence the \(\Delta_{J_{1},J_{2}}^{p}\)'s satisfy the Plücker relations in the homogenous coordinate ring of \(\P(\Wedge^{2}\K^{N_{\delta_{p}}})\), and \(\Delta^{p}\) factors through the Plücker embedding \(\Grass(2,N_{\delta_{p}})\hookrightarrow\P(\Wedge^{2}\K^{N_{\delta_{p}}})\).
Furthermore:
\begin{equation}
  \label{eq:dimSpan}
  B_{p}
  \loceq
  \Set*{
    (\a^{p},\pxi)\in S_{p}\times\pi_{\Omega_{M}}\inverse(V)
    \midbar 
    \dim\Span
    \big(
    {(\mathbi{\alpha}^{p})}_{V}(\a^{p},\xi),
    {(\mathbi{\theta}^{p})}_{V}(\a^{p},\xi)
    \big)
    <2
  }.
\end{equation}

Let us record the following useful information, that follows directly from the construction, see~\eqref{eq:DeltaJJ}.
\begin{lemm}
  \label{lemm:Delta}
  Denote the Plücker line bundle on \(\Grass(2,N_{\delta_{p}})\) by
  \[
    \Q_{\delta_{p}}
    \bydef 
    \O_{\P(\Wedge^{2}\K^{N_{\delta_{p}}})}(1)
    \big\rvert_{\Grass(2,N_{\delta_{p}})}.
  \]
  Outside the indeterminacy locus  of \(\Delta^{p}\), one has
  \[
    (\Delta^{p})^{\ast}
    \Q_{\delta_{p}}
    =
    (\pr_{2}^{[1]})^{\ast}
    \left(
    \O_{\POM}(1)\otimes\pi_{\Omega_{M}}^{\ast}\O_{M}\big(2(\varepsilon_{p}+\delta_{p})\big)
    \right).
  \]
\end{lemm}

The maps \(\Delta^{p}\), \(p=1,\dotsc,c\) will help to define the map \(\varPsi\) mentionned in Sect.~\ref{se:outline}, and thus it is crucial in our proof to control the base locus \(B_{p}\).
The consideration of \eqref {eq:dimSpan} therefore leads us to study the linear maps (\(\pxi\in\pi_{\Omega_{M}}\inverse (V)\))
\[
  \displaymap{\varphi^{p}_{\xi}}{
    S_{p} 
    &\to& 
    \K^{N_{\delta_{p}}}
    \times
    \K^{N_{\delta_{p}}}
    \\
    \a^{p}
    &\mapsto&
    \big(
    {(\mathbi{\alpha}^{p})}_{V}(\a^{p},\xi),
    {(\mathbi{\theta}^{p})}_{V}(\a^{p},\xi)
    \big).
  }
\]
We will need to bound their ranks.
In that aim, we will work with the natural stratification of \(\POM\) induced by contact orders with the \((N+1)\) hypersurfaces \(D_{i}=(\zeta_{i}=0)\).
Recall that the fixed \(\zeta_{0},\dotsc,\zeta_{N}\in H^{0}(M,\O_{M}(1))\) are in general position.
\begin{defi}
  \label{defi:stratification}
  Take \(I\subsetneq\Set{0,\dotsc,N}\) of cardinality \(\#I=N-k_{0}\).
  \begin{enumerate}
    \item 
      We set 
      \(
      D_{I}
      \bydef
      \cap_{i\in I}D_{i}
      \).
      Since \(D_{0},\dotsc,D_{N}\) are in general position,
      \(D_{I}\) is smooth of dimension \(k_{0}\).
    \item 
      We set
      \(
      M_{I}
      \bydef 
      D_{I}\setminus\cup_{i\not\in I}D_{i}
      \).
      One has \(\dim M_{I}=k_{0}\). The different \(M_{I}\) stratify \(M\). 
    \item
      Given  \(I'\subseteq I\subsetneq\Set{0,\dotsc,N}\) of cardinality \(\#I'= N-k_{1}\) we set
      \[
	\Sigma(I,I')
	\bydef
	\P(\Omega_{D_{I'}}\rvert_{M_{I}})
	\setminus 
	\Big(\underset{i\in I\setminus I'}{\textstyle\bigcup}\P(\Omega_{D_{I'\cup\Set{i}}}\rvert_{M_{I}})\Big)
	\subseteq
	\POM.
      \]
      One has \(\dim\Sigma(I,I')=k_{0}+k_{1}-1\). 
      The different \(\Sigma(I,I')\) stratify \(\POM\).
  \end{enumerate}
\end{defi}
Let us justify the last claim of this definition. Observe that locally, 
on a trivializing open set, a point \(\pxi\in\POM\) is in \(\Sigma(I,I')\) if and only if
\[
  \big(
  \zeta_{i}(\xi)=0 
  \iff
  i\in I
  \big)
  \mand
  \big( 
  (\zeta_{i}(\xi)=0, \diff\zeta_{i}(\xi)=0
  \iff
  i\in I'
  \big).
\]
It is thus clear that the different \(\Sigma(I,I')\) stratify \(\POM\).

Given \(I\subsetneq\Set{0,\dotsc,N}\), the inclusion \(\PkI\hookrightarrow\P^{N}\) induces a projection  
\[
  \res_{I}
  \colon
  \K^{N_{\delta_{p}}}
  \cong
  H^{0}\big(\P^{N},\O_{\P^{N}}(\delta_{p})\big)
  \to 
  H^{0}\big(\PkI,\O_{\PkI}(\delta_{p})\big),
\]
identified with the projection
\(
(t_{J})_{\abs{J}=\delta_{p}}
\mapsto
(t_{J})_{\abs{J}=\delta_{p},[J]\cap I=\varnothing}
\), 
\cf~\eqref{eq:Identification}.
\begin{lemm}
  \label{lemm:rk}
  For a point \(\pxi\in\pi_{\Omega_{M}}\inverse(V)\cap\Sigma(I,I')\),
  with \(I'\subseteq I\) of respective cardinalities \(N-k_{1}\leq N-k_{0}\),
  for \(p=1,\dotsc,c\), if \(\delta_{p}\geq 2\),
  \[
    \rk(\varphi_{\xi}^{p})
    \geq
    2\binom{k_{0}+\delta_{p}}{k_{0}}+(k_{0}+1)(k_{1}-k_{0})
    \geq
    2\binom{k_{0}+\delta_{p}}{k_{0}}+(k_{1}-k_{0}).
  \]
  Moreover, \((\res_{I}\times\res_{I})\circ\varphi_{\xi}^{p}\) is a surjective map from \(S_{p}\) to 
  \(H^{0}\big(\PkI,\O_{\PkI}(\delta_{p})\big)\times H^{0}\big(\PkI,\O_{\PkI}(\delta_{p})\big)\).
\end{lemm}
\begin{proof}
  Since \(\alpha_{J}^{p}\) and \(\theta_{J}^{p}\) involve only \(a_{J}^{p}\), 
  the linear map \(\varphi_{\xi}^{p}\) is diagonal by blocks. Each block 
  \(
  H^{0}(M,\O_{M}(\varepsilon_{p}))
  \to
  \K\times\K
  \)
  corresponds to the map
  \(
  a_{J}^{p}
  \mapsto
  \big(
  {(\alpha_{J}^{p})}_{V}(\a^{p},\xi),
  {(\theta_{J}^{p})}_{V}(\a^{p},\xi)
  \big)
  \) 
  for a certain \(J\).
  Recall from~\eqref{eq:theta} that locally (for readability, we now drop the ``local'' notation)
  \[
    \alpha_{J}^{p}(\a^{p},\xi)
    =
    a_{J}^{p}(\xi)(\mathbi{\zeta}(\xi))^{J}
    \mand
    \theta_{J}^{p}(\a^{p},\xi)
    =
    \diff a_{J}(\xi)(\mathbi{\zeta}(\xi))^{J}+(r+1)a_{J}(\xi)\diff(\mathbi{\zeta})^{J}(\xi)
  \]
  There are thus two cases, depending on \(J\): Either \((\mathbi{\zeta}(\xi))^{J}\neq0\), or not.

  Before discussing these two cases, let us note that (since \(\varepsilon_{p}\geq1\) and since sections of very ample line bundles separate first order jets) there exists \(b_{1}\in H^{0}(M,\O_{M}(\varepsilon_{p}))\) such that 
  \(b_{1}(\xi)=0\)
  and
  \(\diff b_{1}(\xi)\neq0\),
  and there also exists \(b_{2}\in H^{0}(M,\O_{M}(\varepsilon_{p}))\) such that 
  \(b_{2}(\xi)\neq0\).

  In the first case, when \((\mathbi{\zeta}(\xi))^{J}\neq0\), the respective images \((0,*_{\neq 0})\) and \((*_{\neq0},*)\) of \(b_{1}\) and \(b_{2}\) cannot be colinear. 
  Thus the rank of the block is \(2\). 
  Observe that there are exactly \(\binom{k_{0}+\delta_{p}}{k_{0}}\) such multi-indices \(J\):
  The assumption \(\pxi\in\Sigma(I,I')\) implies that \(\pi_{\Omega_{M}}(\pxi)\in M_{I}\) and 
  therefore, by definition, \(\mathbi{\zeta}(\xi)^{J}\neq 0\) if and only if \([J]\cap I=\varnothing\). 

  In the other case, when \((\mathbi{\zeta}(\xi))^{J}=0\), clearly 
  \(
  \alpha_{J}^{p}(\a^{p},\xi)
  =
  a_{J}^{p}(\xi)(\mathbi{\zeta}(\xi))^{J}
  =0
  \), whence the rank is at most \(1\).
  We consider only particularly simple instances of this case, for which we can easily compute the rank.
  Since \(\delta_{p}\geq 2\), we can take monomials \(\zeta_{i}^{\delta_{p}-1}\zeta_{j}\) 
  with \(i\in\Set{0,\dotsc,N}\setminus I=\Set{i\midbar\zeta_{i}(\xi)\neq0}\) 
  and 
  \(j\in I\setminus I'=\Set{j\midbar\zeta_{j}(\xi)=0\)
  and
  \(\diff\zeta_{j}(\xi)\neq0}\)% 
  ---it is left to the reader to check that other instances are substantially more complicated---.
  Then the image of \(b_{2}\) is
  \(
  \big(0,(r+1)b_{2}(\xi)\zeta_{i}(\xi)^{\delta_{p}-1}\diff \zeta_{j}(\xi)\big).
  \)
  The entry in the second slot is non zero by assumption on \(i,j\) and \(b_{2}\). Thus the rank of the block is \(1\).

  The announced lower bound for the rank follows, after a count of blocks considered in each case.

  If one consider the composition with the projection onto \(H^{0}\big(\PkI,\O_{\PkI}(\delta_{p})\big)^{\times 2}\), there are only blocks of the first kind, hence the rank is full.
\end{proof}

This lemma implies the following.
\begin{prop} 
  \label{prop:openS}
  As soon as \(\delta_{p}\geq 2\), \(B_{p}\) does not dominate \(S_{p}\) under the canonical projection \(\pr_{1}^{[1]}\). 
  In particular, there exists a non-empty open subset \(S_{p}^{\circ}\subset S_{p}\) such  that the indeterminacy locus of \(\Delta^{p}\) does not intersect \(S_{p}^{\circ}\times\POM\).
\end{prop}
\begin{proof}
  To prove the first assertion, we will establish that for each \(I'\subseteq I\subsetneq\Set{0,\dotsc,N}\) the  locus 
  \[
    Z_{p}(I,I')
    \bydef 
    B_{p}\cap(S_{p}\times\Sigma(I,I'))
    \subset
    S_{p}\times\Sigma(I,I')
  \] 
  does not dominate \(S_{p}\). 
  This will imply that \(B_{p}\) does not dominate \(S_{p}\) under \(\pr_{1}^{[1]}\).
  It suffices then to take \(S_{p}^{\circ}\bydef S_{p}\setminus\pr_{1}^{[1]}(B_{p})\).
  Accordingly, we fix \(I'\subseteq I\), with respective cardinalities \(N-k_{1}\leq N-k_{0}\).

  One has
  \(
  Z_{p}(I,I')
  =
  Z_{p}^{\alpha}(I,I')
  \sqcup
  (Z_{p}(I,I')\setminus Z_{p}^{\alpha}(I,I'))
  \),
  where
  \[
    Z_{p}^{\alpha}(I,I')
    \bydef
    \Set*{
      (\a^{p},\pxi)\in S_{p}\times \Sigma(I,I')
      \midbar
      \alpha_{J}^{p}(\a^{p},\xi)=0,\forall J
    }.
  \]
  To prove that \(Z_{p}(I,I')\) does not dominate \(S_{p}\), it suffices to prove that neither 
  \(Z_{p}^{\alpha}(I,I')\) nor  \(Z_{p}(I,I')\setminus Z_{p}^{\alpha}(I,I')\) dominates \(S_{p}\).

  We first prove the first part of this statement, namely that \(Z_{p}^{\alpha}(I,I')\) does not dominate  \(S_{p}\). 
  Let us define 
  \[
    Z_{p}^{\alpha}(I)
    \bydef
    \Set*{
      (\a^{p},x)\in S_{p}\times M_{I}
      \midbar
      \alpha_{J}^{p}(\a^{p},x)=0,\forall J
    }
    =
    (\id_{S_{p}}\times\pi_{\Omega_{M}})
    (Z_{p}^{\alpha}(I,I')) .
  \] 
  Observe that 
  \(\pr_{1}^{[1]}(Z_{p}^{\alpha}(I,I'))=\pr_{1}(Z_{p}^{\alpha}(I))\),
  thus \(Z_{p}^{\alpha}(I,I')\) dominates \(S_{p}\) if and only if \(Z_{p}^{\alpha}(I)\) does.

  Fix \(x\in M_{I}\). Recall that  \(\dim M_{I}=k_{0}\).
  For any \(\a^{p}=(a_{J}^{p})_{\abs{J}=\delta_{p}}\), one has 
  \[
    \tag{$\ast$}\label{eq:Equivalence}
    \big(
    \alpha_{J}(\a^{p},x)=0,\ 
    \forall J \ \text{such that}\ 
    \abs{J}=\delta_{p}
    \big)
    \iff
    \big(
    a_{J}^{p}(x)=0,\ 
    \forall J \ \text{such that}\ 
    [J]\cap I=\varnothing
    \big). 
  \]
  Indeed, by definition one has \(\alpha_{J}(\a^{p},x)=a_{J}^{p}(x)\mathbi{\zeta}(x)^{J}\). 
  Therefore \(\alpha_{J}(\a^{p},x)=0\) if and only if \(a_{J}^{p}(x)=0\) or \(\mathbi{\zeta}(x)^{J}=0\). 
  But on the other hand, since \(x\in M_{I}\), one has that \(\zeta_{i}(x)=0\) if and only if \(i\in I\). 
  As a consequence we get that \(\zeta(x)^{J}=0\) if and only if \([J]\cap I \neq\varnothing\). 
  This yields the equivalence in~\eqref{eq:Equivalence}.

  Since there are \(\binom{k_{0}+\delta_{p}}{k_{0}}\) multi-indices \(J\) satisfying the right hand side condition in~\eqref{eq:Equivalence}, 
  and since each such multi-index \(J\) imposes an independent  condition \(a_{J}^{p}(x)=0\), we infer that 
  \[
    \dim(Z_{p}^{\alpha}(I))-\dim(S_{p})
    \leq 
    \dim M_{I}+\max_{x\in M_{I}}\dim(Z_{p}^{\alpha}(I)\cap{\pr_{2}\inverse}(x))-\dim(S_{p})
    \leq
    k_{0}-\binom{k_{0}+\delta_{p}}{k_{0}}
    <
    0.
  \]
  Hence, \(Z_{p}^{\alpha}(I)\) does not dominate \(S_{p}\) and neither does \(Z_{p}^{\alpha}(I,I')\).

  We then prove the second part of the sought statement, namely that \(Z_{p}(I,I')\setminus Z_{p}^{\alpha}(I,I')\) does not dominate \(S_{p}\).
  Fix \(\pxi\in \Sigma(I,I')\) and take an open subset \(V\subset M\) containing \(\pi_{\Omega_{M}}(\pxi)\) over which \(\O_{M}(1)\) is trivialized.
  Recall
  \[
    {(\varphi^{p}_{\xi})}_{V}(\a^{p})
    =
    \big(
    {(\mathbi{\alpha}^{p})}_{V}(\a^{p},\xi),
    {(\mathbi{\theta}^{p})}_{V}(\a^{p},\xi)
    \big).
  \]
  Recalling~\eqref{eq:dimSpan}, we have that \((\a^{p},\pxi)\) belongs to  \(Z_{p}(I,I')\) if and only if the dimension of 
  \[
    \Span\big((\mathbi{\alpha}^{p})_{V}(\a^{p},\xi),(\mathbi{\theta}^{p})_{V}(\a^{p},\xi)\big)
    \subset 
    \K^{N_{\delta_{p}}}
  \]
  is strictly less than \(2\).
  Using the stratification, we can even be more precise.
  Since \((\a^{p},\pxi)\not\in Z_{p}^{\alpha}(I,I')\), \ie \({(\mathbi{\alpha}^{p})}_{V}(\a^{p},\xi)\neq0\), this rank condition means that  \[
    (\mathbi{\theta}^{p})_{V}(\a^{p},\xi)
    \in
    \K
    \cdot
    (\mathbi{\alpha}^{p})_{V}(\a^{p},\xi).
  \]
  In particular, \(\alpha_{J}^{p}(\a^{p},\xi)=0\implies\theta_{J}^{p}(\a^{p},\xi)=0\). 

  For shortness, denote \(Z_{\xi}\bydef (\pr^{[1]}_{2})\inverse(\pxi)\cap (Z_{p}(I,I')\setminus Z_{p}^{\alpha}(I,I'))\).
  Remark that \(Z_{\xi}\) is isomorphic to \(\pr^{[1]}_{1}(Z_\xi)\) via \(\pr_{1}^{[1]}\).
  We have just seen that
  \[
    \varphi_{\xi}^{p}\big(\pr_{1}^{[1]}(Z_{\xi})\big)
    \subseteq
    \Set*{
      (v,w)
      \in 
      \K^{N_{\delta_{p}}}\times \K^{N_{\delta_{p}}}
      \midbar 
      w\in \K \cdot v\textbf{ and }v_{J}=w_{J}=0\text{ if }[J]\cap I\neq\varnothing
    }.
  \]
  The dimension of the space on the right hand side being \(\binom{k_{0}+\delta_{p}}{k_{0}}+1\), and the dimension of the kernel of \(\varphi_{\xi}^{p}\) being \(\bigl(\dim S_{p}-\rk(\varphi_{\xi}^{p})\bigr)\), it follows that
  \[
    \dim Z_{\xi}
    =
    \dim \pr^{[1]}_{1}(Z_\xi)
    \leq
    \binom{k_{0}+\delta_{p}}{k_{0}}+1+\dim S_{p}-\rk(\varphi_{\xi}^{p}).
  \]
  Since this argument holds for any \(\pxi\in \Sigma(I,I')\), we obtain 
  \[
    \dim\big(Z_{p}(I,I')\setminus Z_{p}^{\alpha}(I,I')\big)
    \leq
    \binom{k_{0}+\delta_{p}}{k_{0}}+1
    +
    \dim S_{p}
    -
    \min_{\xi\in\Sigma(I,I')}\rk(\varphi_{\xi}^{p})
    +
    \dim(\Sigma(I,I')).
  \]
  Taking account of the different dimensions and of the above Lemma~\ref{lemm:rk}, this yields that
  \[
    \dim\big(Z_{p}(I,I')\setminus Z_{p}^{\alpha}(I,I')\big)
    -
    \dim S_{p}
    \leq
    \binom{k_{0}+\delta_{p}}{k_{0}}
    +
    1
    -
    2\binom{k_{0}+\delta_{p}}{k_{0}}
    -
    (k_{1}-k_{0})
    +
    k_{0}+k_{1}-1
    \leq
    2k_{0}-\binom{k_{0}+\delta_{p}}{k_{0}}.
  \]
  It is left to the reader to check that since \(\delta_{p}\geq 2\), the last expression is negative. 
  Therefore \(Z_{p}(I,I')\setminus Z_{p}^{\alpha}(I,I')\) does not dominate \(S_{p}\) and this finishes the proof.
\end{proof}

For each \(1\leq p\leq c\), when \(\delta_{p}\geq2\), by Proposition~\ref{prop:openS} the restriction of \(\Delta^{p}\) to \(S_{p}^{\circ}\times\POM\) is a morphism and not merely a rational map.
Let us briefly recall from~\eqref{eq:defDelta} that
\(\Delta^{p}(\a^{p},\pxi)\) corresponds locally to the \(2\)-plane in \(H^{0}(\P^{N},\O_{\P^{N}}(\delta_{p}))\) 
spanned by the polynomials 
\(
\sum_{\abs{J}=\delta_{p}}
\alpha_{J}^{p}(\a^{p},\xi)\mathbi{z}^{J}
\)
and
\(
\sum_{\abs{J}=\delta_{p}}
\theta_{J}^{p}(\a^{p},\xi)\mathbi{z}^{J}
\).
Proposition~\ref{prop:openS} ensures that if \(\a^{p}\in S_{p}^{\circ}\), those two polynomials are indeed linearly independent.
Assume \(\delta_{1},\dotsc,\delta_{c}\geq 2\) and denote
\[
  S^{\circ}
  \bydef 
  \big(S_{1}^{\circ}\times\dotsb\times S_{c}^{\circ}\big)
  \cap
  S.
\]
We get a map
\begin{equation}
  \label{eq:DefinitionPsi}
  \displaymap{\varPsi}{
    S^{\circ}\times\POM
    &\to& 
    \Grass(2,N_{\delta_{1}})
    &\times\dotsb\times&
    \Grass(2,N_{\delta_{c}})
    &\times&
    \P^{N}
    \\
    (\a,\pxi)
    &\mapsto&
    \big( \Delta^{1}(\a^{1},\pxi)
    &,\dotsc,&
    \Delta^{c}(\a^{c},\pxi)
    &,&
    [(\mathbi{\zeta}(\xi))^{r}]\big),
  }
\end{equation}
where 
\([(\mathbi{\zeta}(\xi))^{r}]\)
stands for 
\([\zeta_{0}(\xi)^{r}:\dotso:\zeta_{N}(\xi)^{r}]\).

For the sake of shortness, denote
\[
  \G
  \bydef
  \Grass(2,N_{\delta_{1}})
  \times
  \dotsb
  \times
  \Grass(2,N_{\delta_{c}}),
\]
and consider the universal family
\begin{equation}
  \label{eq:Y}
  \Y
  \bydef
  \Set*{
    (\Delta_{1},\dotsc, \Delta_{c},z)
    \in 
    \G
    \times 
    \P^{N}
    \midbar\\
    \forall P\in \Delta_{1}\cup\dotsb\cup\Delta_{c}
    \colon 
    P(z)=0
  }.
\end{equation}
Along 
\(
\P(\Omega_{\X/S^{\circ}})
\bydef
\POXS
\cap
\big(S^{\circ}\times\POM\big)
\),
the map \(\varPsi\) factors through \(\Y\). 
Indeed, above an affine open subset \(V\subset M\) with a fixed trivialization of \(\O_{M}(1)\rvert_{V}\),
the variety \(\P(\Omega_{\X/S^{\circ}})\cap(S^{\circ}\times\pi_{\Omega_{M}}\inverse(V))\) is defined by the \(c\) couples of equations 
\((E^{p})_{V}=\sum (\alpha^{p}_{J})_{V}({\mathbi{\zeta}})_{V}^{rJ}\) 
and 
\((\diff E^{p})_{V}=\sum(\theta^{p}_{J})_{V}({\mathbi{\zeta}})_{V}^{rJ}\), 
for \(p=1,\dotsc,c\). 
A straightforward computation shows that 
\(\varPsi\big(\P(\Omega_{\X/S^{\circ}})\cap(S^{\circ}\times\pi_{\Omega_{M}}\inverse(V))\big)\subset\Y\).
This can be summarized in the following commutative diagram:
\begin{equation}
  \label{eq:Diagram}
  \vcenter{\xymatrix{
    S^{\circ}\times\POM
    \ar[r]^-{\varPsi} 
    \ar@{<-_{)}}[d(0.75)]
    &\G\times\P^{N}
    \ar[r]^<<<{\rho}
    \ar@{<-_{)}}[d(0.75)]
    &\G
    \\
    \P(\Omega_{\X/S^{\circ}})
    \ar[r]^-{\varPsi}
    &\Y\ar[ru]^{\rho}
  }}.
\end{equation}
Here we denote by \(\rho\colon\G\times\P^{N}\to\G\) the projection onto the first factor, 
as well as its restriction to \(\Y\subset\G\times\P^{N}\).

Observe that, since \(2c\geq N\), the projection \(\rho\colon\Y\to\G\) is generically finite, and might be not surjective.
One can however think of \(\Y\), at least formally, as the universal family of ``complete intersections'' of codimension \(2c\) and multi-degree \((\delta_{1},\delta_{1},\dotsc,\delta_{c},\delta_{c})\) parametrized by \(\G\).
Having this in mind, thanks to a result of Benoist (see Sect.~\ref{se:Olivier} below),  we will obtain a key ingredient needed in our proof to pull back the positivity by \(\varPsi\), namely the fact that the locus of elements in \(\G\) that parametrize  positive dimensional  schemes is of large codimension in \(\G\).

For technical reasons, we will now consider the stratification of \(\POM\)
induced by the vanishing of the sections \(\zeta_{0},\dotsc,\zeta_{N}\) on the base \(M\). 
Namely, for \(M_{I}\) defined as above (see Defintion~\ref{defi:stratification}), the different
\[
  \POMI
  =
  \pi_{\Omega_{M}}\inverse(M_{I})
\]
stratify \(\POM\).
For reasons that will soon become obvious, we are mainly interested in the case where \(\dim(M_{I})=k\geq1\) \ie in sets \(I\) of cardinality \(N-k\leq N-1\). 
Fix \(I\subsetneq\Set{0,\dotsc,N}\) with cardinality \(\#I<N\) and take \((\a,\pxi)\in\P(\Omega_{\X/S^{\circ}})\cap(S^{\circ}\times\POMI)\). 
By definition, as \(\pxi\in\POMI=\pi_{\Omega_{M}}\inverse(M_{I})\) we obtain immediately that \((\mathbi{\zeta}(\xi))^{r}\in\PkI\subset\P^{N}\) 
(recall that \(\K_{I}=\Set{(z_{0},\dotsc,z_{N})\in \K^{N+1}\midbar z_{i}=0,\forall i\in I}\)).
In other words, along
\(
\P(\Omega_{\X/S^{\circ}})\cap(S^{\circ}\times\POMI)
\), 
the map 
\(\varPsi\)
factors through
\[
  \Y(I)
  \bydef
  \Set*{
    (\Delta_{1},\dotsc, \Delta_{c},z)\in 
    \G
    \times 
    \PkI
    \midbar
    \forall P\in \Delta_{1}\cup\dotsb\cup\Delta_{c}
    \colon 
    P(z)=0
  }
  =
  \Y\cap(\G\times\PkI).
\]
Let us denote by \(\rho_{I}\colon\Y(I)\to\G\) the restriction to \(\Y(I)\) of the canonical projection \(\G\times\PkI\to\G\).

\subsection{Avoiding the exceptional locus}
\label{sse:exc}
Because we want to use the positivity of the base \(\G\), it is important to avoid positive dimensional fibers of \(\rho_{I}\colon\Y(I)\to \G\). We denote the \textsl{exceptional locus} of a generically finite morphism \(f\colon X\to Y\) between two varieties by
\[
  \Exc(f)
  \bydef
  \Set*{
    x\in X
    \midbar
    \dim_{x}\big(f\inverse(\Set{f(x)})\big)>0
  }.
\]
Note that for any \(I\subsetneq\Set{0,\dotsc,N}\), \(\Exc(\rho_{I})\subseteq\Y(I)\subseteq\G\times\P^{N}\).

\begin{lemm}
  \label{lemm:HyperplaneRestriction}
  If \(\delta_{1},\dotsc,\delta_{c}\geq\dim(\POM)\),
  then there exists a non-empty open subset \(U\subset S^{\circ}\) such that
  for each \(I\subsetneq\Set{0,\dotsc,N}\) with cardinality \(\#I<N\),
  \[
    \big(U\times\POMI\big)
    \cap 
    \varPsi\inverse(\Exc(\rho_{I}))
    =
    \varnothing.
  \]
\end{lemm}
\begin{proof}
  It is sufficient to prove that for each \(I\subsetneq\Set{0,\dotsc,N}\) with cardinality \(N-k\), for \(k\geq1\), 
  there exists a non-empty open subset \(U(I)\subset S^{\circ}\) such that 
  \[
    \big(U(I)\times\POMI\big)
    \cap 
    \varPsi\inverse(\Exc(\rho_{I}))
    =
    \varnothing.
  \]
  We will prove the sought statement by dimension count, establishing that 
  \[    
    \dim\Big(
    \big(S^{\circ}\times\POMI\big)\cap{\varPsi}\inverse\big(\Exc (\rho_{I})\big)
    \Big)
    <
    \dim S
    =
    \dim S^{\circ}.
  \]

  Let us fix a subset \(I\) as above and set
  \(
  \G^{\infty}
  \bydef
  \rho_{I}\big(\Exc(\rho_{I})\big)
  \subseteq
  \G
  \).
  Since \(\rho_{I}\) is the restriction of \(\rho\) to \(\Y(I)\subseteq\Y\), one has 
  \(
  \Exc(\rho_{I})
  \subseteq
  \rho\inverse(\G^{\infty})
  \) 
  and therefore
  \[
    \varPsi\inverse(\Exc(\rho_{I}))
    \subseteq
    \varPhi\inverse(\G^{\infty}),
  \]
  where (\cf~\eqref{eq:DefinitionPsi}),
  \[
    \varPhi
    \bydef
    \rho\circ\varPsi
    =
    \Delta^{1}\times\dotsb\times\Delta^{c}
    \colon
    S^{\circ}\times\POM
    \to 
    \G.
  \]
  Fix \(\pxi\in\POMI\) and denote by \(\varPhi_{\xi}\colon S^{\circ}\to\G\) the map defined by \(\varPhi_{\xi}(\a)\bydef\varPhi(\a,\pxi)\).
  Observe that the projection \(\pr_{1}^{[1]}\) induces an isomorphism between 
  \((S^{\circ}\times\Set{\pxi})\cap\varPhi\inverse(\G^{\infty})\) and \(\varPhi_{\xi}\inverse(\G^{\infty})\).
  We will prove the inequality
  \[
    \dim(\varPhi_{\xi}\inverse(\G^{\infty}))<\dim S-\dim\POMI,
  \]
  which implies the lemma, since then
  \[
    \dim S
    >
    \dim\Big(
    \big(S^{\circ}\times\POMI\big)\cap{\varPhi}\inverse(\G^{\infty})
    \Big)
    \geq
    \dim\Big(
    \big(S^{\circ}\times\POMI\big)\cap{\varPsi}\inverse(\Exc(\rho_{I}))
    \Big).
  \]

  In order to study \(\dim(\varPhi_{\xi}\inverse(\G^{\infty}))\), let us first define an affine analogue to the family \(\rho_{I}\colon\Y(I)\to\G\).
  Denote
  \[
    \widetilde{\G}
    \bydef
    H^{0}(\PkI,\O_{\PkI}(\delta_{1}))^{\times 2}
    \times
    \dotsb
    \times
    H^{0}(\PkI,\O_{\PkI}(\delta_{c}))^{\times 2}.
  \]
  An affine analogue of \(\Y(I)\) is
  \[
    \widetilde{\Y}(I)
    \bydef
    \Set*{
      (\alpha_{1},\theta_{1},\dotsc,\alpha_{c},\theta_{c},z)
      \in 
      \widetilde{\G}\times\PkI
      \midbar
      \alpha_{i}(z)=0,\theta_{i}(z)=0,\forall i=1,\dotsc,c
    }.
  \]
  Let us denote by \(\tilde{\rho}_{I}\colon\widetilde{\Y}(I)\to\widetilde{\G}\) the restriction to \(\widetilde{\Y}(I)\) of the canonical projection \(\widetilde{\G}\times\PkI\to\widetilde{\G}\).

  Next, we define a natural analogue to \(\varPhi_{\xi}\).
  Fix an open chart \(V\subset M\) containing \(\pi_{\Omega_{M}}(\pxi)\) and trivializing \(\O_{M}(1)\). 
  Then, with the notation introduced before Definition~\ref{defi:stratification}, 
  consider the  linear maps 
  \(
  \tilde{\varphi}_{\xi}^{p}
  \bydef
  (\res_{I}\times\res_{I})\circ
  \varphi_{\xi}^{p}
  \),
  for \(1\leq p\leq c\) .
  These induce a linear map 
  \(
  \widetilde{\varPhi}_{\xi}
  \bydef 
  \tilde{\varphi}_{\xi}^{1}\times\dotsb\times\tilde{\varphi}_{\xi}^{c}
  \)
  from \(S\) to \(\widetilde{\G}\),
  which is surjective by Lemma~\ref{lemm:rk}.

  We are now ready to conclude.
  Define
  \(
  \widetilde{\G}^{\infty}
  \bydef
  \tilde{\rho}_{I}\big(\Exc(\tilde{\rho}_{I})\big)\subset\widetilde{\G}
  \).
  Observe that for any \(\a\in S^{\circ}\), one has
  \[
    \rho_{I}\inverse(\varPhi_{\xi}(\a))
    =
    \tilde{\rho}_{I}\inverse(\widetilde{\varPhi}_{\xi}(\a))
    \subset
    \PkI.
  \] 
  From this it follows that 
  \(\varPhi_{\xi}(\a)\in \G^{\infty}\)
  if and only if 
  \(\widetilde{\varPhi}_{\xi}(\a)\in \widetilde{\G}^{\infty}\).
  Hence
  \[
    \varPhi_{\xi}\inverse(\G^{\infty})
    = 
    \widetilde{\varPhi}_{\xi}\inverse(\widetilde{\G}^{\infty})\cap S^{\circ}.
  \]
  But now, since \(\widetilde{\varPhi}_{\xi}\) is linear and surjective, it has equidimensional fibers of dimension \(\dim S-\dim \widetilde{\G}\).
  Thus
  \[
    \dim\big(\widetilde{\varPhi}_{\xi}\inverse(\widetilde{\G}^{\infty})\big)
    =
    \dim S-\dim\widetilde{\G}+\dim\widetilde{\G}^{\infty}
    =
    \dim S-\codim(\widetilde{\G}^{\infty},\widetilde{\G}).
  \]
  Besides, one has 
  \[
    \codim(\widetilde{\G}^{\infty},\widetilde{\G})
    >
    \min_{1\leq p\leq c}(\delta_{p})
    \geq 
    \dim\POM
    \geq
    \dim\POMI.
  \]
  The first  inequality follows from Corollary~\ref{coro:CodimBound}, to be proven below.
  In the notation of Sect.~\ref{se:Olivier}, apply the result to \(k=2c\geq N\) and \((e_{1},\dotsc,e_{k})=(\delta_{1},\delta_{1},\dotsc,\delta_{c},\delta_{c})\). In this case \(\widetilde{\G}=A_{\mathbi{e}}\) and \(\widetilde{\G}^{\infty}=A_{\mathbi{e}}^{\infty}\), thus we have  
  \(
  \codim(\widetilde{\G}^{\infty},\widetilde{\G})
  >
  \min_{1\leq p\leq c}(\delta_{p})
  \).
  The second inequality follows from the hypothesis on \(\delta_{p}\).
  From the above considerations, one deduces that 
  \[
    \dim\big(\varPhi_{\xi}\inverse(\G^{\infty})\big)
    \leq
    \dim\big(\widetilde{\varPhi}_{\xi}\inverse(\hat{\G}^{\infty})\big)
    <
    \dim S
    -
    \POMI.
    \qedhere
  \]
\end{proof}

\subsection{Pulling back the positivity}
\label{sse:positivity}
We are now in position to implement the last part of our strategy, thereby establishing the following stronger version of Theorem~\ref{thm:ample}.
\begin{theo}
  \label{thm:ampleStrong}
  For each choice of \(\delta_{1},\dotsc,\delta_{c}\geq\dim\POM\), there exists an integer \(m(\mathbi{\delta})\) such that 
  for any
  \(r> 2m(\mathbi{\delta})(\abs{\mathbi{\varepsilon}}+\abs{\mathbi{\delta}})\), a general member of the family \(\X=\X(\mathbi{\delta},\mathbi{\varepsilon},r)\) constructed in Sect.~\ref{sse:bihom} has ample cotangent bundle.
\end{theo}
By the openness property of ampleness, this immediately implies the following.
\begin{coro}
  \label{coro:ampleDague}
  On a N-dimensional smooth projective variety \(M\), equipped with a very ample line bundle \(\O_{M}(1)\), if \(N/2\leq c\leq N\), for any degrees \((d_{1},\dotsc,d_{c})\in(\Zpos)^{c}\) satisfying
  \begin{equation}
    \label{eq:condition}
    \tag{\ensuremath{\dagger}}
    {}^{\exists} 
    \delta_{1},\dotsc,\delta_{c}\geq2N-1,
    {\ }^{\exists} 
    \varepsilon_{1},\dotsc,\varepsilon_{c}\geq1,
    {\ }^{\exists}
    r>
    2m(\mathbi{\delta})(\abs{\mathbi{\varepsilon}}+\abs{\mathbi{\delta}})
    \colon
    \quad
    d_{p}
    =
    \delta_{p}
    (r+1)
    +
    \varepsilon_{p}\quad
    {\scriptstyle(p=1,\dotsc,c)},
  \end{equation}
  the complete intersection 
  \(
  X
  \bydef
  H_{1}\cap\dotsb\cap H_{c}
  \)
  of general hypersurfaces 
  \(
  H_{1}\in\abs*{\O_{M}(d_{1})},\dotsc,H_{c}\in\abs*{\O_{M}(d_{c})}
  \)
  has ample cotangent bundle.
\end{coro}
Condition~\eqref{eq:condition} will be discussed in Sect.~\ref{sse:condition} below.
\begin{proof}[Proof of Theorem~\ref{thm:ampleStrong}]
  For \(\a\in S\), denote \(\pi_{\a}\colon\P(\Omega_{X_{\a}})\to X_{\a}\). 
  It is sufficient to prove that for a sufficiently large \(\mu\),
  and \(\a\) general,
  the line bundle
  \(\O_{\P(\Omega_{X_{\a}}(-1/\mu))}(\mu)=\O_{\P(\Omega_{X_{\a}})}(\mu)\otimes\pi_{\a}^{\ast}\O_{X_{\a}}(-1)\) 
  is nef.
  For such \(\mu\) and \(\a\), we will thus prove that any irreducible curve \(\mathscr{C}\subset\P(\Omega_{X_{\a}})\) satisfies
  \[
    \label{eq:NefGoal}
    \tag{\ensuremath{\ast}}
    \mathscr{C}
    \cdot 
    \left(
    \O_{\P(\Omega_{X_{\a}})}(\mu)
    \otimes
    \pi_{\a}^{\ast}\O_{X_{\a}}(-1)
    \right)
    \geq
    0.
  \]

  Observe that, since \(\big(\POMI\big)_{I\subsetneq\Set{0,\dotsc,N}}\) stratifies  \(\POM\), for each irreducible curve \(\mathscr{C}\subset \P(\Omega_{X_{\a}})\subset \POM\), there exists a unique set \(I\) such that \(\POMI\) contains an open dense subset of \(\mathscr{C}\).
  We fix \(\a\in S\), such a curve \(\mathscr{C}\) and the corresponding \(I\), of cardinality \(N-k\).
  If \(\dim(M_{I})=k=0\), \eqref{eq:NefGoal} holds for all \(\mu\), so we assume that \(k\geq 1\).

  We fix \(\mathbi{\delta}=(\delta_{1},\dotsc,\delta_{c})\) and we work with the ample line bundle
  \[
    \Q
    \bydef
    \Q_{\delta_{1}}
    \boxtimes 
    \dotsb
    \boxtimes
    \Q_{\delta_{c}}
    \longrightarrow
    \G
    =
    \Grass(2,N_{\delta_{1}})
    \times
    \dotsb
    \times
    \Grass(2,N_{\delta_{c}}),
  \]
  where for \(\delta\in\N\),  we denote by \(\Q_{\delta}\) the Pl\"{u}cker line bundle on \(\Grass(2,N_{\delta})\).
  Let us denote by \(\rho_{I}\) the restriction of \(\rho\) to \(\Y(I)\).
  Since \(\Q\) is ample and \(\rho_{I}\colon\Y(I)\to\G\) is generically finite, the pullback-bundle \({\rho_{I}}^{\ast}\Q\) is big and nef. Accordingly, its \textsl{null locus}---which, thanks to a theorem due to Nakamaye, is  known to coincide with the \textsl{augmented base locus} \(\B_{+}({\rho_{I}}^{\ast}\Q)\), see~\cite{Nak00} or \eg~\cite{Laz04II} Sect.~10.3---is exactly the reunion of all positive dimensional fibers of \(\rho_{I}\), \ie the exceptional locus of the map \(\Y_{I}\to\G\).

  From the very definition of \(\B_{+}\), the fact that \(\Q\boxtimes \O_{\PkI}(1)\) is very ample, and Noetherianity, we know that there exists an integer \(m_{I}(\mathbi{\delta})\) such that for \(m\geq m_{I}(\mathbi{\delta})\), 
  \[
    \B_{+}(\rho_{I}^{\ast}\Q)
    =
    \Bs\left(({\Q}^{\otimes m}\boxtimes \O_{\PkI}(-1))\rvert_{\Y(I)}\right).
  \]

  To sum up, there exists \(m_{I}(\mathbi{\delta})\) such that for \(m\geq m_{I}(\mathbi{\delta})\)
  \[
    \Exc(\rho_{I})
    =
    \Bs\left(({\Q}^{\otimes m}\boxtimes \O_{\PkI}(-1))\rvert_{\Y(I)}\right)
    =
    \Bs\left(({\Q}^{\otimes m}\boxtimes \O_{\P^{N}}(-1))\rvert_{\Y(I)}\right).
  \]

  We will now fix \(\mu\), depending only on \(\mathbi{\delta}\).
  In order to work uniformly on \(\a\) and \(\mathscr{C}\), set 
  \[
    m(\mathbi{\delta})
    \bydef
    \max\Set*{
      m_{I}(\mathbi{\delta})
      \midbar 
      I\subsetneq\Set{0,\dotsc,N},\#I\leq N-1
    }.
  \]
  Recall from Lemma~\ref{lemm:Delta} that for \(p=1,\dotsc,c\):
  \[
    (\Delta^{p}\rvert_{S^{\circ}})^{\ast}
    \Q_{\delta_{p}}
    =
    (\pr_{2}^{[1]})^{\ast}
    \left(\O_{\POM}(1)\otimes\pi_{\Omega_{M}}^{\ast}\O_{M}\big(2(\varepsilon_{p}+\delta_{p})\big)\right).
  \]
  Using the definition~\eqref{eq:DefinitionPsi} of \(\varPsi\), restricted to \(\P(\Omega_{\X/S^{\circ}})\), this yields
  \begin{equation}
    \label{eq:PullBack}
    \varPsi^{\ast}
    \left(
    \Q^{\otimes m(\mathbi{\delta})}\boxtimes\O_{\P^{N}}(-1)
    \right)
    =
    \O_{\P(\Omega_{\X/S^{\circ}})}(cm(\mathbi{\delta}))
    \otimes
    \pi_{\Omega_{\X/S^{\circ}}}^{\ast}
    \O_{\X}\Big(-r+2m(\mathbi{\delta})(\abs{\mathbi{\varepsilon}}+\abs{\mathbi{\delta}})\Big),
  \end{equation}
  Observe the \(-r\); our assumption on \(r\) is precisely \(-r+2m(\mathbi{\delta})(\abs{\mathbi{\varepsilon}}+\abs{\mathbi{\delta}})<0\).

  We take \(\mu\bydef cm(\mathbi{\delta})\).
  Then, we proceed as follows.
  Since by assumption \(\delta_{1},\dotsc,\delta_{c}\geq\dim(\POM)\), the conclusion of Lemma~\ref{lemm:HyperplaneRestriction} is satisfied for any \(r\). 
  Hence for \(\a\in U\), 
  for every irreducible curve \(\mathscr{C}\subset\P(\Omega_{X_{\a}})\)
  having an open dense part in \(\Set{\a}\times\POMI\)
  with \(\#I<N\), one has necessarily
  \[
    \mathscr{C}
    \not\subseteq 
    \varPsi\inverse(\Exc(\rho_{I}))
    =
    \varPsi\inverse\big(\Bs\big((\Q^{\otimes m(\mathbi{\delta})}\boxtimes \O_{\PkI}(-1))\rvert_{\Y(I)}\big)\big).
  \]

  This implies the existence of a section 
  \(\sigma\in H^{0}(\Y(I),\Q^{\otimes m(\mathbi{\delta})}\boxtimes\O_{\PkI}(-1))\)
  such that \(\mathscr{C}\not\subseteq(\varPsi^{\ast}\sigma=0)\). 
  In particular,
  \(\mathscr{C}\cdot (\varPsi^{\ast}\sigma=0)\geq 0\).
  Therefore by~\eqref{eq:PullBack}, 
  \[
    \mathscr{C}
    \cdot 
    \left(
    \O_{\P(\Omega_{X_{\a}})}(\mu)
    \otimes
    \pi_{\a}^{\ast}
    \O_{X_{\a}}\left(-r+m(\mathbi{\delta})(\abs{\mathbi{\varepsilon}}+\abs{\mathbi{\delta}})
    \right)
    \right)
    \geq
    0.
  \]
  This proves~\eqref{eq:NefGoal} and thus finishes the proof.
\end{proof}

\subsection{The condition on the degrees}
\label{sse:condition}
Up to taking a multiple, one can assume that all \(\delta_{i}\) in Theorem~\ref{thm:ample} are larger than \(2N-1\). Then by Corollary~\ref{coro:ampleDague}, the degrees \(d_{i}\) for which the conclusion of  Theorem~\ref{thm:ample} holds in the direction \(\mathbi{\delta}\), are the ones satisfying the following condition. 
\[
  \tag{\ensuremath{\dagger}}
  {}^{\exists} 
  \varepsilon_{1},\dotsc,\varepsilon_{c}\geq1,
  {\ }^{\exists}
  r>
  2m(\mathbi{\delta})(\abs{\mathbi{\varepsilon}}+\abs{\mathbi{\delta}})
  \colon
  \quad
  d_{p}
  =
  \delta_{p}
  (r+1)
  +
  \varepsilon_{p}\quad
  {\scriptstyle(p=1,\dotsc,c)}.
\]

The picture below shows particular solutions, with \(1\leq\varepsilon_{p}\leq\delta_{p}\), of~\eqref{eq:condition} in an arbitrary direction \(\mathbi{\delta}=(\delta_{1},\dotsc,\delta_{c})\), with \(\delta_{p}\geq2N-1\).
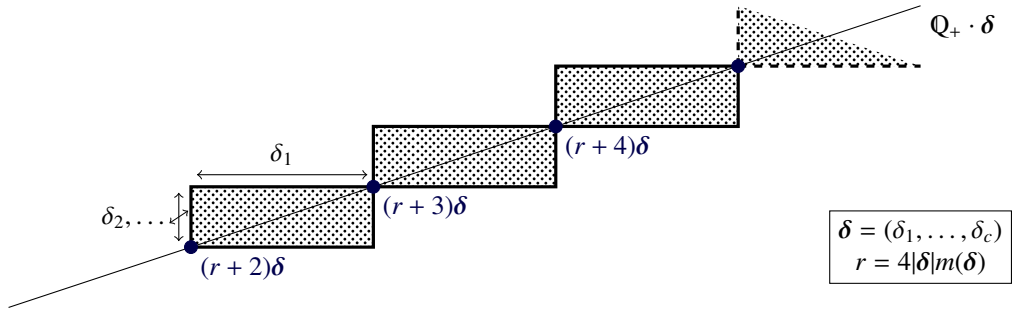
\begin{figure}[htbp]
  \centering
  \begin{tikzpicture}[scale=.8]
    \draw (0,0)--(15,5) node[below right]{$\mathbb{Q}_{+}\cdot\mathbi{\delta}$};
    \draw[very thick,pattern=crosshatch dots] (3,1) rectangle (6,2);
    \draw[<->] (3.1,2.2)--(5.9,2.2) node[midway, above]{$\delta_{1}$};
    \draw[<->] (2.8,1.1)--(2.8,1.9) node[midway, left]{$\delta_{2},\dotsc$};
    \draw[<->,very thin] (2.65,1.4)--(2.95,1.6);
    \draw[very thick,pattern=crosshatch dots] (6,2) rectangle (9,3);
    \draw[very thick,pattern=crosshatch dots] (9,3) rectangle (12,4);
    \draw[very thick, dashed] (12,4)--(15,4) (12,4)--(12,5);
    \fill[pattern=crosshatch dots] (12,4)--(12,5)--(15,4)--cycle;
    \filldraw[blue!30!black] (3,1) circle(3pt) node[below right]{$(r+2)\mathbi{\delta}$};
    \filldraw[blue!30!black] (6,2) circle(3pt)node[below right]{$(r+3)\mathbi{\delta}$};
    \filldraw[blue!30!black] (9,3) circle(3pt)node[below right]{$(r+4)\mathbi{\delta}$};
    \filldraw[blue!30!black] (12,4) circle(3pt);
    \node[draw, thin, align=center]at (15,1){$\mathbi{\delta}=(\delta_{1},\dotsc,\delta_{c})$\\$r=4\abs{\mathbi{\delta}}m(\mathbi{\delta})$};
  \end{tikzpicture}
  \caption{Some degrees satisfying the condition~\eqref{eq:condition}}
\end{figure}
\newline \noindent
As one can see, these solutions lie in a union of rectangular blocks that contains the ray spanned by \(\mathbi{\delta}\) starting from a sufficiently large multiple \(\nu(\mathbi{\delta})\cdot\mathbi{\delta}\) with \eg \(\nu(\mathbi{\delta})=4\abs{\mathbi{\delta}}m(\mathbi{\delta})+2\).
So, in particular, we treat all large degrees \((d_{1},\dotsc,d_{c})\) in the ray spanned by \((\delta_{1},\dotsc,\delta_{c})\), with a lower bound depending on the direction.

\section{On the codimension of the exceptional locus}
\label{se:Olivier}
Given \(N\geq 1\), if one parametrizes complete intersections in \(\P^{N}\) by products of spaces of homogeneous polynomials, the locus of families of polynomials not parametrizing complete intersections is ``small'' compared to the entire parameter space. This follows easily from the work of Benoist (\cite{Ben11}), but since it is playing a central role in our proof, we provide all the details here.
We do not claim any originality in this section. 
\begin{lemm}[{\cite[Lemme 2.3]{Ben11}}]
  Let \(X\subset\P^{N}\) be a (irreducible) subvariety of dimension \(n\). Let \(\mathcal{G}\) be the set of all hypersurfaces of degree \(e\) containing \(X\). Then 
  \[
    \codim(\mathcal{G},\abs*{\O_{\P^{N}}(e)})
    \geq 
    \binom{e+n}{n}.
  \]
\end{lemm}
\begin{proof}
  This is the proof from~\cite{Ben11}, translated in English for the reader's convenience. 

  Let \(L\) be a \((N-n-1)\)-dimensional linear subspace of \(\P^{N}\) such that \(L\cap X=\varnothing\). 
  Let \(\pi_{L}\colon X\to\P^{n}\) be the projection induced by \(L\). 
  All the fibers of \(\pi_{L}\) are non-empty and finite. This proves that if \(\mathcal{C}\) is the set of all cones of degree \(e\) with vertex \(L\), \(\mathcal{C}\cap \mathcal{G}=\varnothing\), and therefore
  \[
    \codim(\mathcal{G},\abs*{\O_{\P^{N}}(e)})
    \geq
    \dim\mathcal{C}+1=\binom{e+n}{n}.
    \qedhere
  \]
\end{proof}

For \(k\geq 1\) and integers \(e_{1},\dotsc,e_{k}\) denote
\(
T_{\mathbi{e}}
\bydef 
\abs*{\O_{\P^{N}}(e_{1})}
\times
\dotsb
\times
\abs*{\O_{\P^{N}}(e_{k})}
\), 
and define \(T_{\mathbi{e}}^{\nci}\) to be the subvariety of \(T_{\mathbi{e}}\) parametrizing equations that do not define a complete intersection variety if \(k \leq N\), or to be the subvariety parametrizing equations that define a positive dimensional scheme if \(k\geq N\). 
Observe that for \(k\leq N\), the locus \(T_{\mathbi{e}}^{\nci}\) is also the locus parametrizing equations that define a subscheme of \(\P^{N}\) whose dimension exceeds the expected dimension \(N-k\). 
Thus in the case \(k=N\) the two notions coincide. 
In any case \(T_{\mathbi{e}}^{\nci}\) is closed in \(T_{\mathbi{e}}\). 
\begin{coro}
  \label{coro:Benoist}
  With the above notation,
  \[
    \codim(T_{\mathbi{e}}^{\nci},{T_{\mathbi{e}}})
    \geq 
    \min_{1\leq i\leq \min\Set{N,k}}\binom{e_{i}+N-i+1}{N-i+1}
    >
    \min_{1\leq i\leq\min\Set{N,k}}(e_{i}).
  \]
  In particular, this codimension tends to infinity as the \(e_{i}\)'s tend to infinity.
\end{coro}
\begin{proof}
  The proof is an induction on \(k\). If \(k=1\), it follows from the fact that any hypersurface in \(\P^{N}\) is a complete intersection. Suppose the result holds for \(k<N\). Then for an element \(F_{1},\dotsc,F_{k+1}\in T_{\mathbi{e}}\), the subscheme \(X_{k+1}\) defined by \(F_{1},\dotsc,F_{k+1}\) fails to be a complete intersection, either when the subscheme \(X_{k}\) defined by the first \(k\) equations fails to be a complete intersection, or when the hypersurface \((F_{k+1}=0)\) contains one of the irreducible components of \(X_{k}\). The codimension of the subset defined  by the first condition is greater than \(\min_{1\leq i\leq k}\binom{e_{i}+N-i+1}{N-i+1}\) by induction hypothesis. The codimension of the subset defined by the second condition is greater than \( \binom{e_{k+1}+N-k}{N-k}\). This is seen simply by fixing some \(X_{k}\) which is a complete intersection, so that all irreducible components of \(X_{k}\) are of dimension \(N-k\) and applying Benoist's Lemma to each of those irreducible components. This proves the statement for all \(1\leq k\leq N\), the statement for \(k>N\) follows at once. Indeed, if \(k>N\) then the codimension of \(T^{\nci}_{\mathbi{e}}\) in \(T_{\mathbi{e}}\) is larger than the codimension of the set of elements in \(T_{\mathbi{e}}\) whose  \(N\) first equations define a positive dimensional subscheme of \(\P^{N}\), and we are reduced to the statement for \(k=N\).
\end{proof}

Let us lastly observe that Corollary~\ref{coro:Benoist} allows us to deduce a similar bound when one considers the affine analogue of \(T_{\mathbi{e}}\) given by 
\[
  A_{\mathbi{e}}
  \bydef
  H^{0}\big(\P^{N},\O_{\P^{N}}(e_{1})\big)
  \times
  \dotsb
  \times
  H^{0}\big(\P^{N},\O_{\P^{N}}(e_{k})\big).
\] 
We state the result in the case \(k\geq N\geq1\), which is sufficient for our application.
Define 
\[
  A_{\mathbi{e}}^{\infty}
  \bydef 
  \Set*{
    (F_{1},\dotsc,F_{k})
    \in 
    A_{\mathbi{e}}
    \midbar 
    \dim\big((F_{1}=0)\cap\dotsb\cap(F_{k}=0)\big)>0
  },
\]
\begin{coro}
  \label{coro:CodimBound}
  With the above notation,
  \[
    \codim(A_{\mathbi{e}}^{\nci},{A_{\mathbi{e}}})
    >
    \min_{1\leq i\leq k}(e_{i}).
  \]
\end{coro}
\begin{proof}
  Denote in this proof
  \(
  A_{\mathbi{e}}^{\circ}
  \bydef
  \prod_{i=1}^{k}
  \big(H^{0}(\P^{N},\O_{\P^{N}}(e_{i}))\setminus\Set{0}\big)
  \). 
  Observe that
  \[
    \codim(A_{\mathbi{e}}^{\infty},A_{\mathbi{e}})
    \geq 
    \min\Set*{
      \codim(A_{\mathbi{e}}^{\infty}\cap A_{\mathbi{e}}^{\circ},A_{\mathbi{e}}^{\circ}),
      \codim (A_{\mathbi{e}}\setminus A_{\mathbi{e}}^{\circ},A_{\mathbi{e}})
    }.
  \]

  Now, on the one hand, from Corollary~\ref{coro:Benoist},
  \[
    \codim(A_{\mathbi{e}}^{\infty}\cap A_{\mathbi{e}}^{\circ},A_{\mathbi{e}}^{\circ})
    =
    \codim(T_{\mathbi{e}}^{\infty},T_{\mathbi{e}})
    > 
    \min_{1\leq i\leq N}(e_{i})
    \geq
    \min_{1\leq i\leq k}(e_{i}).
  \]
  Indeed the canonical projection 
  \(p\colon A_{\mathbi{e}}^{\circ}\to T_{\mathbi{e}}\)
  has equidimensional fibers of dimension \(k\) and   
  \(p\inverse(T_{\mathbi{e}}^{\infty})=A_{\mathbi{e}}^{\infty}\cap A_{\mathbi{e}}^{\circ}\).
  On the other hand, by definition, one has
  \[
    \codim (A_{\mathbi{e}}\setminus A_{\mathbi{e}}^{\circ},A_{\mathbi{e}})
    =
    \min_{1\leq i\leq k} \dim H^{0}\big(\P^{N},\O_{\P^{N}}(e_{i})\big)
    =
    \min_{1\leq i\leq k}\binom{e_{i}+N}{N}
    >
    \min_{1\leq i\leq k}(e_{i}).
    \qedhere
  \]
\end{proof}

\backmatter
\paragraph*{Acknowledgments} 
We would like to thank Olivier Benoist for the conversations we had about the content of this paper, and for allowing us to reproduce some of his arguments. 
We would like to also thank Simone Diverio and Gianluca Pacienza for their friendly support and for useful suggestions.
We gratefully appreciated further useful suggestions of Joël Merker and of Christophe Mourougane.

\bibliographystyle{smfalpha}
\bibliography{BD}
\vfill
\end{document}